\RequirePackage[l2tabu, orthodox]{nag}
\RequirePackage{snapshot}

\documentclass[10pt,twocolumn]{article}

\makeatletter
\if@twocolumn
  \usepackage[dvips,pdftex,letterpaper,top=0.5in, bottom=0.5in, left=0.75in, right=0.5in,includefoot,heightrounded]{geometry}
\else
  \usepackage[dvips,pdftexletterpaper,margin=1in,includefoot,heightrounded]{geometry}
\fi

\usepackage{srcltx}

\usepackage{amsmath}
\usepackage{amssymb,amsfonts}

\usepackage{abstract}

\usepackage[xdvi]{graphicx}
\usepackage[usenames,dvipsnames]{color}
\usepackage[usenames,dvipsnames]{xcolor}
\usepackage{subfigure}

\usepackage{booktabs}

\usepackage{setspace}
\usepackage{flushend}
\usepackage{multicol}

\usepackage{cite}
\usepackage{url}
\usepackage[normalem]{ulem}

\usepackage{enumerate}

\usepackage{hyperref}

\usepackage[sc,tiny]{titlesec}
\usepackage[mathscr]{euscript}

       \usepackage{mathptmx}

\newtheorem{proposition}{Proposition}

\def\QED{\mbox{$\square$}}
\def\proof{\noindent{\it Proof:~}}
\def\endproof{\hspace*{\fill}~\QED\par\endtrivlist\unskip}



\title{%
The Fourier-Like and Hartley-Like Wavelet Analysis Based on Hilbert Transforms
}

\author{%
L. R. Soares%
\thanks{%
L. R. Soares
was
with
Graduate Program in Electrical Engineering,
Universidade Federal de Pernambuco.
}
\and
H. M. de Oliveira%
\thanks{
H. M. de Oliveira 
was
with the Departamento de Eletr\^onica e Sistemas,
Universidade Federal de Pernambuco.
Now he is with
the Signal Processing Group,
Departamento de Estat\'{\i}stica, 
Universidade Federal de Pernambuco.
E-mail: \url{hmo@de.ufpe.br}}
\and
R.~J.~Cintra%
\thanks{%
R.~J.~Cintra
was
with
Graduate Program in Electrical Engineering,
Universidade Federal de Pernambuco.
He is currently
with 
the Signal Processing Group,
Departamento de Estat\'{\i}stica, 
Universidade Federal de Pernambuco.
E-mail: 
\protect\url{rjdsc@de.ufpe.br}
}
}

\date{}

\newcommand{\myabstract}{%
In continuous-time wavelet analysis, most wavelet present some kind of symmetry. Based on the Fourier and Hartley transform kernels, a new wavelet multiresolution analysis is proposed. This approach is based on a pair of orthogonal wavelet functions and is named as the Fourier-Like and Hartley-Like wavelet analysis. A Hilbert transform analysis on the wavelet theory is also included.
}

\newcommand{\mykeywords}{%
Wavelet design, Hilbert transform, Fourier transform kernel, Hartley transform kernel, Fourier-Like wavelets, Analytic wavelets, Hartley-Like wavelets, Continuous-time, Wavelet transform.
}

\begin{document}

\makeatletter
\if@twocolumn

\twocolumn[%
  \maketitle
  \begin{onecolabstract}
    \myabstract
  \end{onecolabstract}
  \begin{center}
    \small
    \textbf{Keywords}
    \\\medskip
    \mykeywords
  \end{center}
  \bigskip
]
\saythanks

\else

  \maketitle
  \begin{abstract}
    \myabstract
  \end{abstract}
  \begin{center}
    \small
    \textbf{Keywords}
    \\\medskip
    \mykeywords
  \end{center}
  \bigskip
  \onehalfspacing
\fi

\section{Introduction}

The idea of comparing Fourier analysis with wavelet decompositions is the starting point for introducing an analysis based on a couple of orthogonal wavelet functions, one with even symmetry and other with odd symmetry. This approach is here presented and named as the Fourier-Like and Hartley-Like wavelet analysis.
In the standard Fourier analysis, a signal $x(t)$ is simultaneously analyzed by even and odd (in quadrature) functions, being represented by:
\begin{align*}
x(t)
=
\text{d.c. term}
+
\text{cosine terms}
+
\text{sine terms}
.
\end{align*}
In the standard wavelet multiresolution analysis (WMRA)~\cite{ref1}, $x(t)$ may be represented by:
\begin{align*}
x(t) = \text{$\phi$ term}
+
\text{$\psi$ terms}
,
\end{align*}
where 
``$\phi$~term'' account for the analysis of $x(t)$ with a scaling function $\phi(t)$
and
the ``$\psi$~terms'' represent those ones derived from scaled versions of a mother wavelet function $\psi(t)$.

Comparing the WMRA and the Fourier analysis equations, the scaling coefficients, ``$\phi$~term”, play a role that corresponds to the d.c. term of the Fourier series; the wavelet coefficients, ``$\psi$~terms'', can be viewed as the harmonic components of the Fourier series, since the harmonics are scaled versions of the infinite Fourier kernel.

Hence, it seems natural to associate to each wavelet with even symmetry another one of odd symmetry, and vice-versa. The Hilbert transform can naturally be invoked to derive the in quadrature version of a symmetric (or an anti-symmetric) wavelet. Hence, x(t) may be represented on a new continuous-time WMRA by:

\begin{align*}
x(t) = 
\text{$\phi$~term}
+ 
\text{$\psi$~terms}
+
\text{orthogonal of $\psi$~terms}
.
\end{align*}

Accordingly, new wavelet functions that look like the Fourier and Hartley transform kernels are invoked. So, by analogy to Fourier and Hartley transforms, the “cosine and sine” kernel is replaced by “$\psi$ and the Hilbert transform of $\psi$” in this new concept of wavelet analysis.
In order to allow further investigation on the Fourier-Like and Hartley-Like wavelet analysis, a brief review of the Hilbert transform is presented, as well the results of applying some wavelet properties to the Hilbert transform of wavelets.

\section{The Hilbert Transform}

The Hilbert transform of a function $g(t)$ is defined by~\cite{ref3}:
\begin{align*}
\operatorname{\mathcal{H}}
\left\{
g
\right\}
(t)
=
\operatorname{p.v.}
\frac{1}{\pi}
\int_{-\infty}^\infty
\frac{g(\tau)}
{t  - \tau}
\operatorname{d}\tau
,
\end{align*}
where 
$\operatorname{p.v.}$
is the Cauchy principal value of the integral. 
After a change of variable, the Hilbert transform can be written as a convolution:
\begin{align}
\label{equation-5}
\operatorname{\mathcal{H}}
\left\{
g
\right\}
(t)
=
\frac{1}{\pi}
\ast
g(t)
.
\end{align}
By taking the Fourier transform of~\eqref{equation-5}, we have:
\begin{align}
\label{equation-6}
\operatorname{\mathcal{F}}
\left\{
\operatorname{\mathcal{H}}
\left[
g
\right]
\right\}
=
-
j
\cdot
\operatorname{sgn}(\omega)
\cdot
\operatorname{\mathcal{F}}
\left\{
g
\right\}
(\omega)
,
\end{align}
where
$\operatorname{\mathcal{F}}$
is the Fourier transform operator and
$\operatorname{sgn}(\cdot)$ is the signum function defined by:
\begin{align}
\label{equation-7}
\operatorname{sgn}(\omega)
=
\begin{cases}
1, & \omega > 0, \\
0, & \omega = 0, \\
-1,& \omega < 0
.
\end{cases}
\end{align}

It follows from~\eqref{equation-6}--\eqref{equation-7}
 that the Hilbert transform of a function imposes a null at 
$\omega=0$ and a
$-\pi/2$ phase delay on the frequency response of that function. 
Other interesting properties of the Hilbert transform are~\cite{ref3}:
\begin{itemize}
\item
A function and its Hilbert transform are orthogonal over the infinite interval;
\item
The Hilbert transform of a real function is a real function;
\item
The Hilbert transform of an even function is an odd function, and vice-versa.
\end{itemize}
In the framework of wavelets, the Hilbert transform of a symmetric (or anti-symmetric) real wavelet is a real anti-symmetric (or symmetric) function. However, it is necessary to verify whether the resulting function is also a wavelet function or not.

\section{The Hilbert Transform on The Wavelet Analysis}

A function $\psi(t)$ is a mother-wavelet, if and only if,
(i) $\psi(t)$ is in the space of finite energy functions
$L^2(\mathbb{R})$, 
and
(ii) $\psi(t)$ satisfies the admissibility condition~\cite{ref1}.
Some properties are explored in the following propositions in view of applying Hilbert transform to wavelets.
Let
$\Psi(\omega)$ be the Fourier transform of $\psi(t)$, 
$\operatorname{\mathcal{H}}\{\psi\}(t)$ be the Hilbert transform of $\psi(t)$,
$\operatorname{\mathcal{E}}[\psi(t)]$ as the energy of $\psi(t)$, 
and 
$\operatorname{C}[\psi(t)]$
as the admissibility coefficient of $\psi(t)$.

\begin{proposition}
\label{proposition-1}
If 
$\psi(t)$ is a real wavelet, 
then 
$\operatorname{\mathcal{H}}\{\psi\}(t)$
is also a real wavelet with same energy 
and
admissibility coefficient of the generating wavelet~$\psi(t)$.
\end{proposition}
\proof
If $\psi(t)$ is a real wavelet, then $\psi(t)$ belongs to~$L^2(\mathbb{R})$ and satisfies the admissibility condition.
Invoking Parseval's theorem, the energy and admissibility coefficient of 
$\operatorname{\mathcal{H}}\{\psi\}(t)$
are given by:
\begin{align*}
\operatorname{\mathcal{E}}
\left[
\operatorname{\mathcal{H}}\{\psi\}(t)
\right]
&
=
\frac{1}{2\pi}
\int_{-\infty}^\infty
\left|
-j
\cdot
\operatorname{sgn}(\omega)
\cdot
\Psi(\omega)
\right|^2
\operatorname{d}\omega
,
\\
\operatorname{C}
[
\operatorname{\mathcal{H}}\{\psi\}(t)
]
&=
\int_{-\infty}^\infty
\frac{\left|
-j
\cdot
\operatorname{sgn}(\omega)
\cdot
\Psi(\omega)
\right|^2}
{|\omega|}
\operatorname{d}\omega
.
\end{align*}
A simple manipulation yields to
\begin{align*}
\operatorname{\mathcal{E}}
\left[
\operatorname{\mathcal{H}}\{\psi\}(t)
\right]
&
=
\frac{1}{2\pi}
\int_{-\infty}^\infty
\left|
\Psi(\omega)
\right|^2
\operatorname{d}\omega
,
\\
\operatorname{C}
[
\operatorname{\mathcal{H}}\{\psi\}(t)
]
&=
\int_{-\infty}^\infty
\frac{\left|
\Psi(\omega)
\right|^2}
{|\omega|}
\operatorname{d}\omega
.
\end{align*}
Applying Parseval's theorem on the right-side of the energy equation, it is straightforward to conclude that 
$\operatorname{\mathcal{H}}\{\psi\}(t)$ 
is also in 
$L^2(\mathbb{R})$. 
Moreover, $\psi(t)$ and $\operatorname{\mathcal{H}}\{\psi\}(t)$ have the same energy.
As $\psi(t) \in L^2(\mathbb{R})$ 
and 
$|\Psi(0)|=0$,
it promptly follows that $\operatorname{\mathcal{H}}\{\psi\}(t)$ also satisfies the admissibility condition: 
$\operatorname{C}
[
\operatorname{\mathcal{H}}\{\psi\}(t)
]<+\infty$. 
Furthermore, $\psi(t)$ and $\operatorname{\mathcal{H}}\{\psi\}(t)$
have the same admissibility coefficient.
\endproof

\begin{proposition}
\label{proposition-2}
Let $\psi(t)$ be a wavelet with $N$~vanishing moments, then 
$\operatorname{\mathcal{H}}\{\psi\}(t)$ has at least $N$~vanishing moments.
\end{proposition}
\proof
The $n$th moment of $\psi(t)$ is defined by~\cite{ref1}:
$M_n[\psi(t)] = 
\int_{-\infty}^\infty t^n \cdot \psi(t) \operatorname{d}t$.
As $\psi(t)$ has $N$ vanishing moments, then 
$M_b[\psi(t)]=0$, for $n = 0,1,\ldots,N-1$.
In the frequency domain, 
the moments of $\psi(t)$ are expressed by~\cite{ref3}:
$M_n[\psi(t)]=\frac{\Psi^{(n)}(0)}{(-2\pi j)^n}$,
where the superscript $(n)$
denotes the $n$th derivative of $\Psi(\omega)$.
Hence, the $n$th moment of $\operatorname{\mathcal{H}}\{\psi\}(t)$ is given by
$$
M_n
\left[
\operatorname{\mathcal{H}}\{\psi\}(t)
\right]
=
\left.
\frac{[-j \cdot \operatorname{sgn}(\omega) \cdot \Psi(\omega)]^{(n)}}{(2\pi j)^n}
\right|_{\omega=0}
.
$$
Consequently, 
it follows that $\operatorname{\mathcal{H}}\{\psi\}(t)$ has at least 
$N$~vanishing moments.
\endproof

\begin{figure*}
\centering

\subfigure[]{\includegraphics[width=0.3\linewidth]{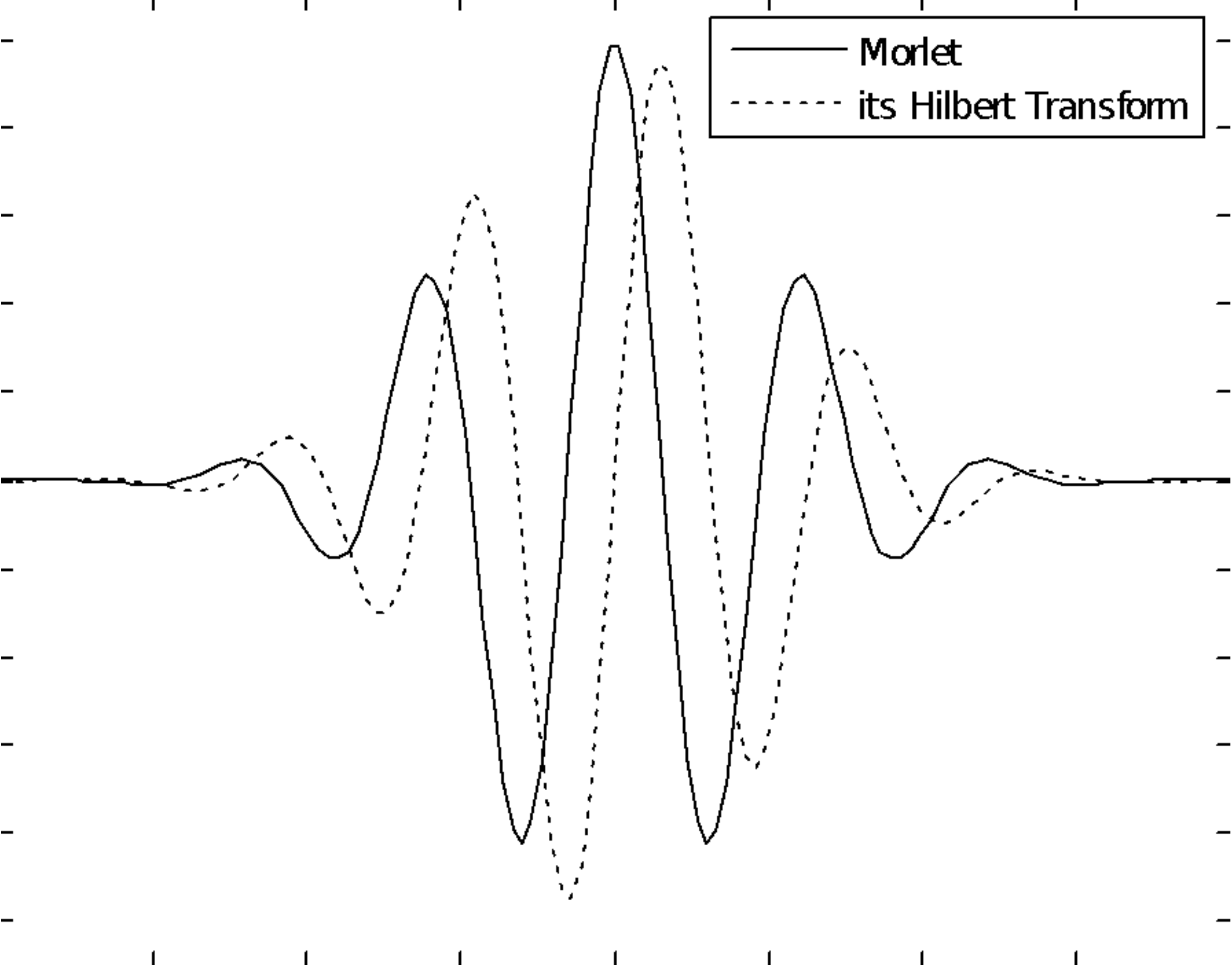}}
\quad
\subfigure[]{\includegraphics[width=0.3\linewidth]{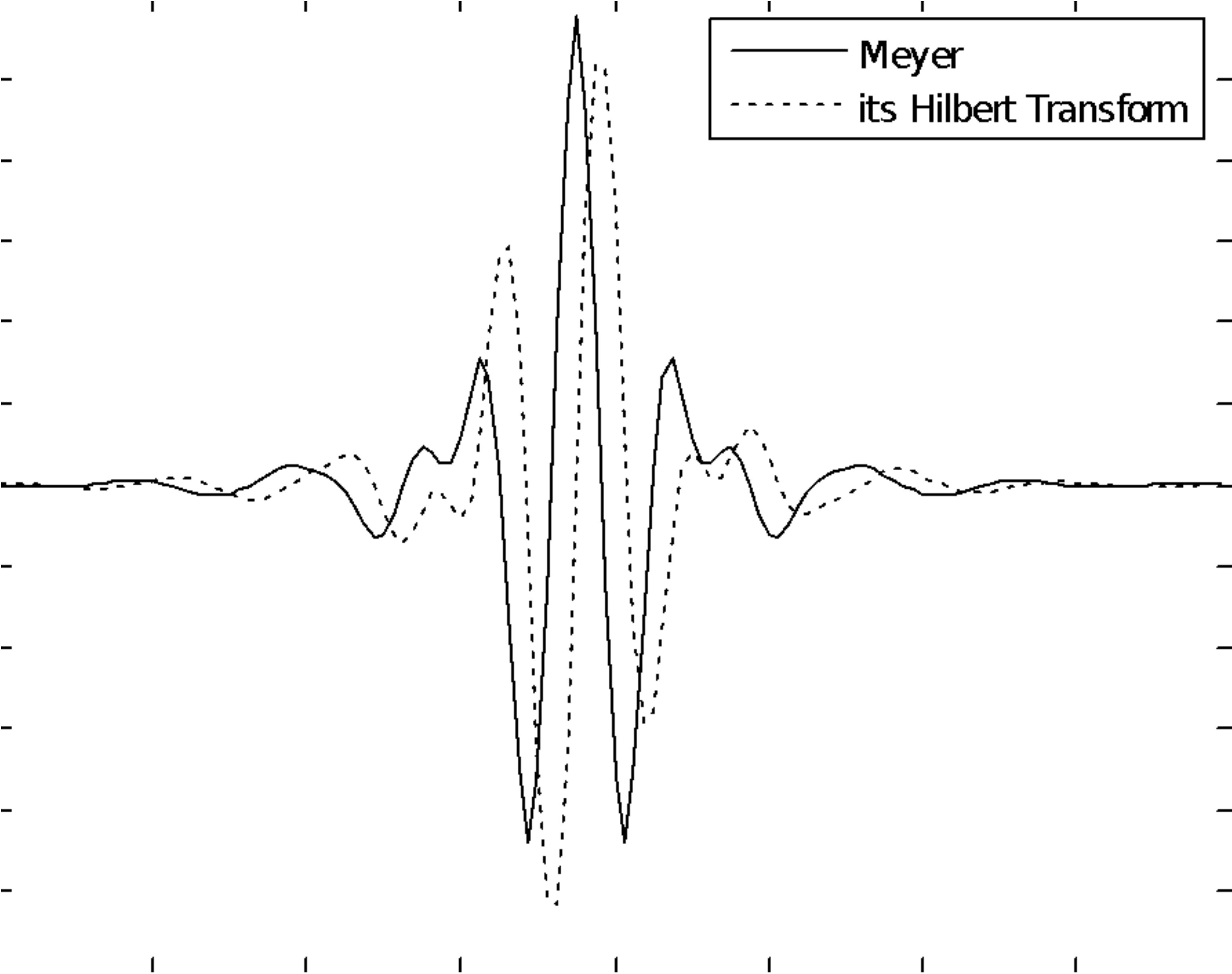}}
\quad
\subfigure[]{\includegraphics[width=0.3\linewidth]{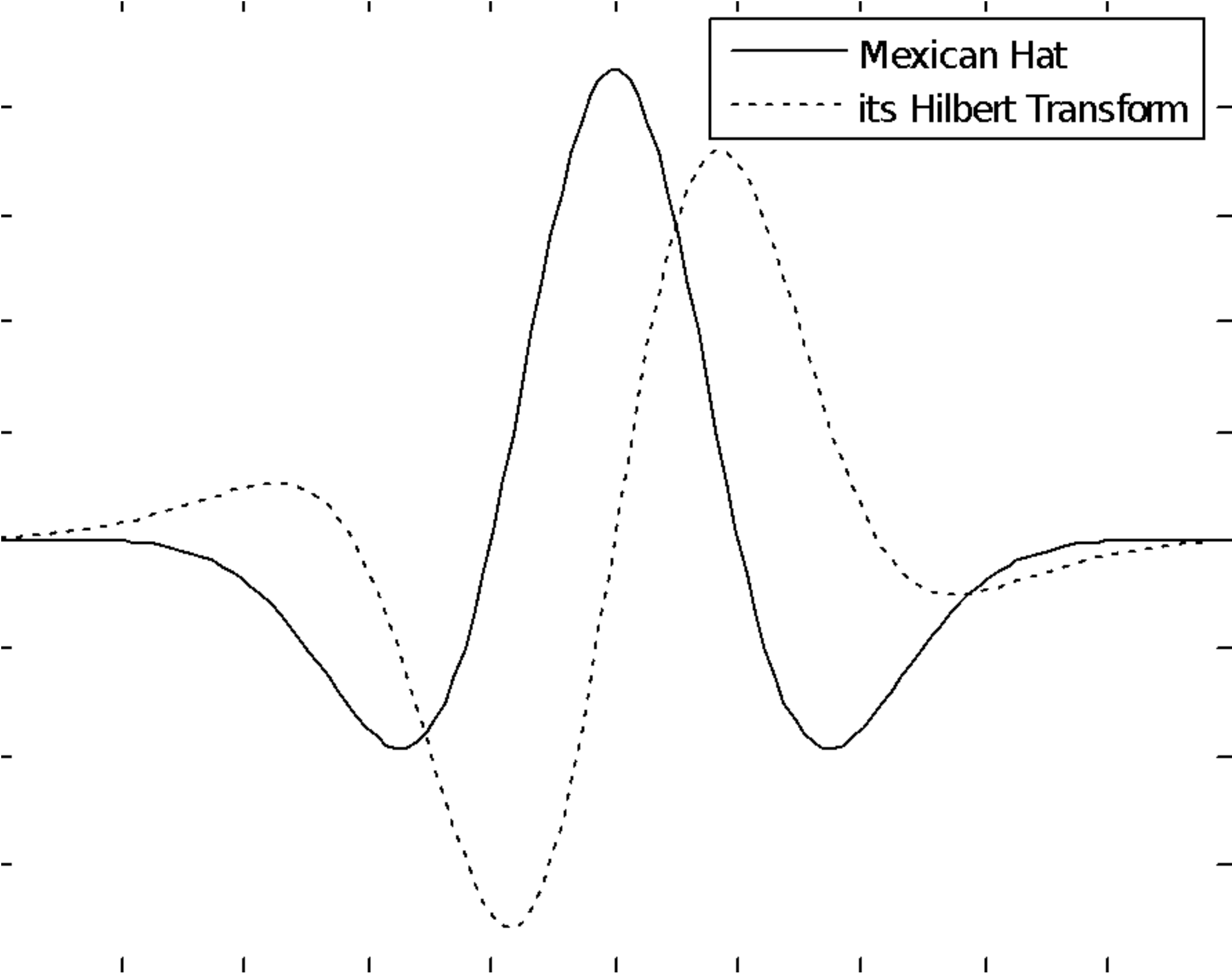}}
\\
\subfigure[]{\includegraphics[width=0.3\linewidth]{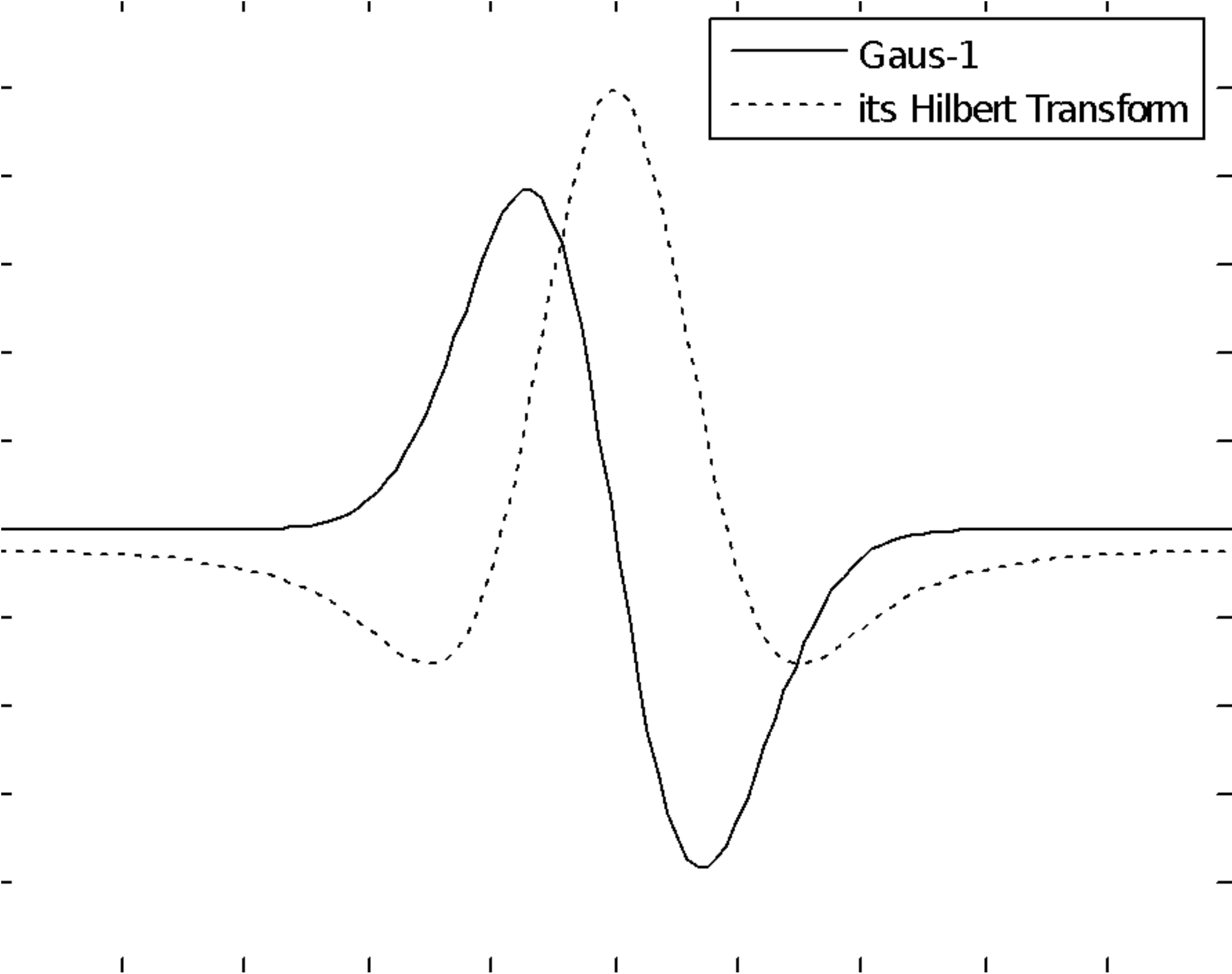}}
\quad
\subfigure[]{\includegraphics[width=0.3\linewidth]{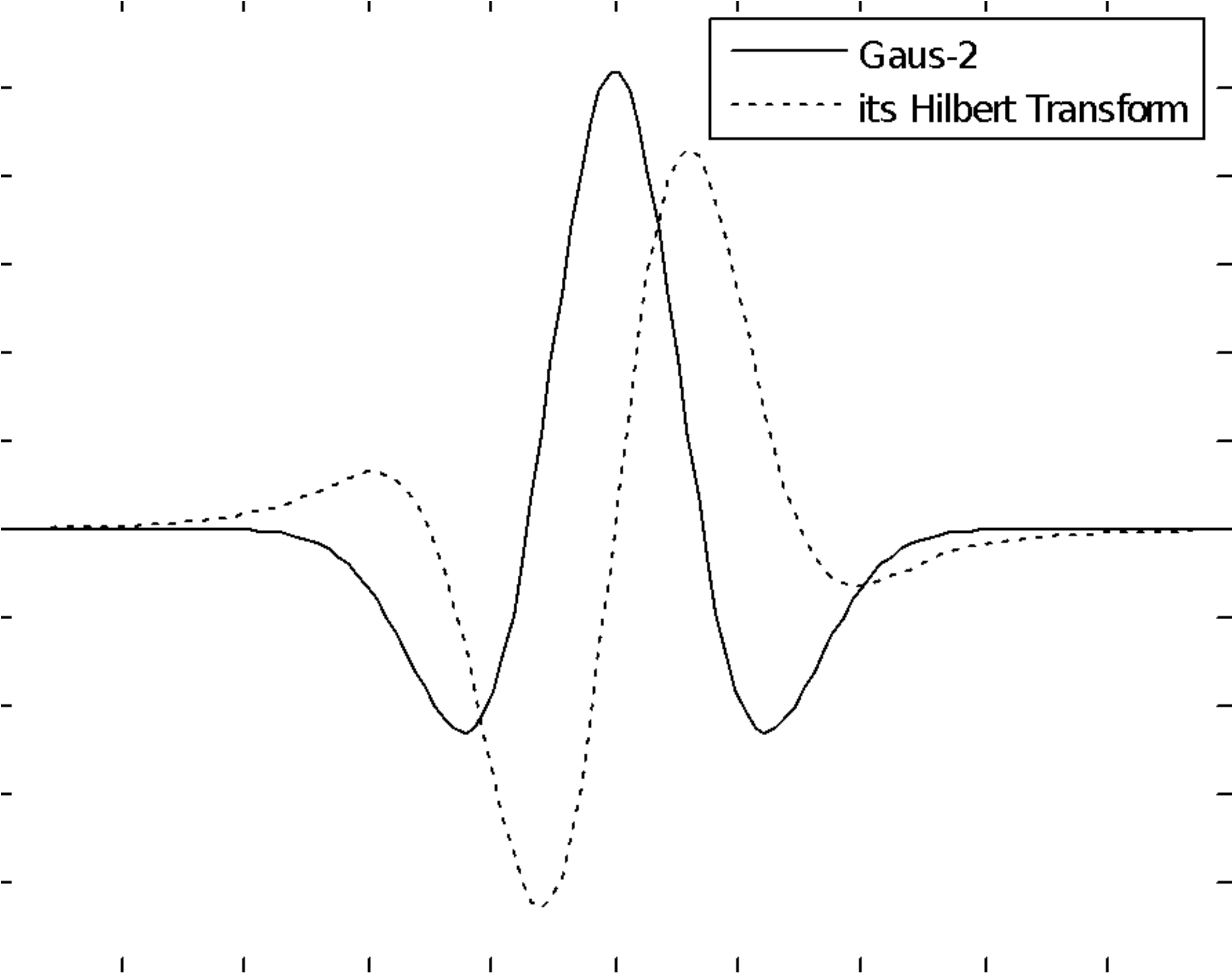}}
\quad
\subfigure[]{\includegraphics[width=0.3\linewidth]{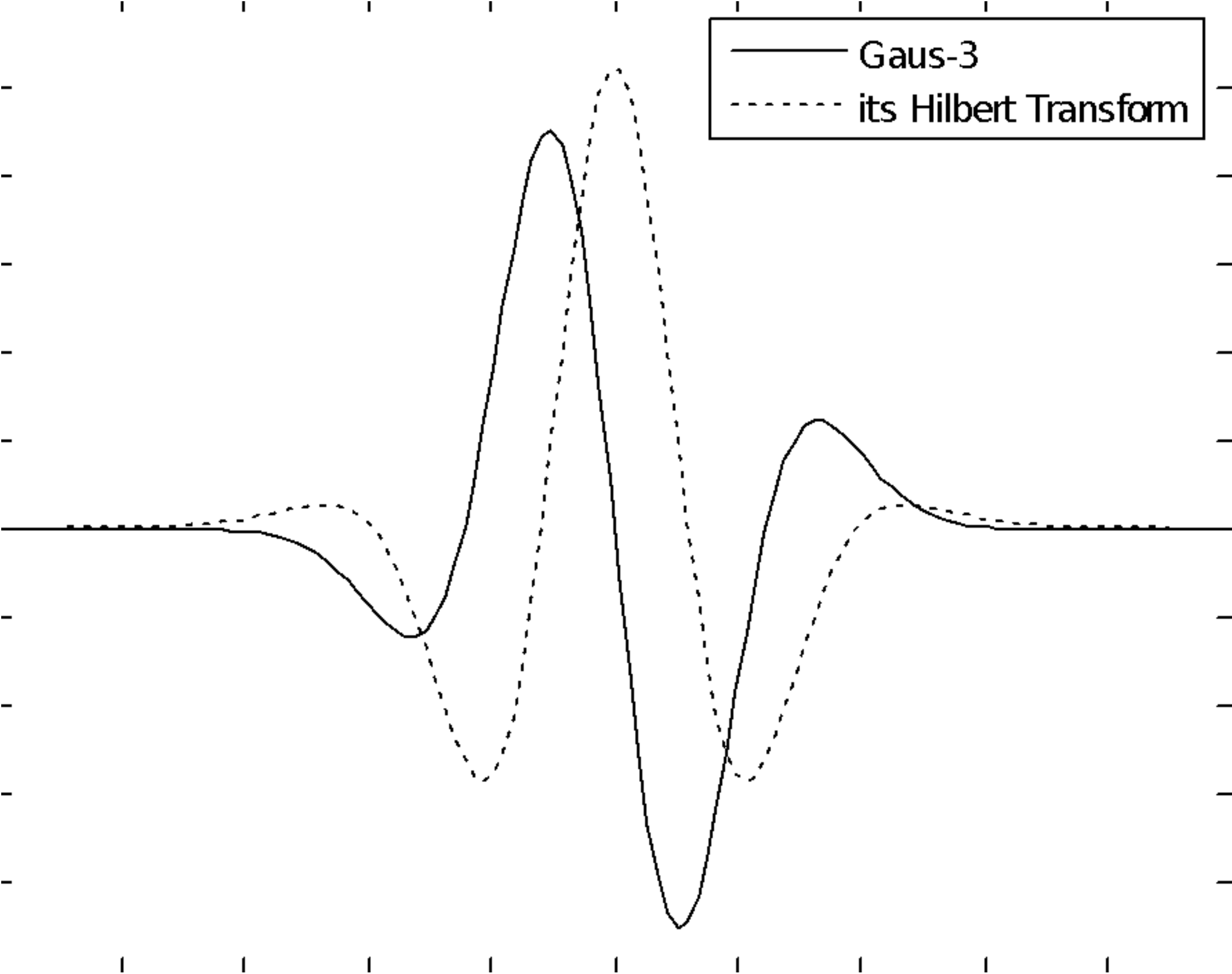}}

\caption{Continuous-time wavelets and their Hilbert transforms: 
(a)~Morlet;
(b)~Meyer;
(c)~Mexican Hat;
(d)~Gaussian-1;
(e)~Gaussian-2;
and
(f)~Gaussian-3.
Scale of the horizontal axis is
one~unit (time) per division.
Scale of the vertical axis is
$0.2$ units (amplitude) per division,
ranging from $-1$ to~$1$.}
\label{figure-1}
\end{figure*}

In view of Propositions~\ref{proposition-1} and~\ref{proposition-2},
it follows that $\operatorname{\mathcal{H}}\{\psi\}(t)$ 
is a wavelet with same energy, admissibility coefficient and at
least same number of null moments than 
its generating wavelet $\psi(t)$.
Figure~\ref{figure-1}
shows a few continuous-time real wavelets and their corresponding Hilbert transforms.

\section{The Fourier Kernel on The Wavelet Analysis}

It is time to define a new wavelet function that looks like the Fourier transform kernel, which can analyze both symmetries of an asymmetric signal.
The Fourier transform kernel, or Fourier kernel, is defined by
$e^{jt} = \cos(t) + j\cdot\sin(t)$.
Considering that
$\operatorname{\mathcal{H}}\{\cos\}(t)=-\sin(t)$,
then the Fourier kernel can also be written as 
$e^{jt} = \cos(t) - j \operatorname{\mathcal{H}}\{\cos\}(t)$.

This naive observation motivates the definition of Fourier-Like wavelets based on a real wavelet and its Hilbert transform. Let us define the Fourier-Like wavelet, 
$\operatorname{\mathscr{F}}\{\psi(t)\}$, by:
\begin{align*}
\operatorname{\mathscr{F}}\{\psi(t)\}
=
\frac{1}{\sqrt2}
\cdot
\left[
\psi(t)
-
j
\cdot
\operatorname{\mathcal{H}}\{\psi\}(t)
\right]
.
\end{align*}

The following proposition
demonstrates that 
$\operatorname{\mathscr{F}}\{\psi(t)\}$ is also a wavelet and that the factor $1/\sqrt{2}$
is imposed so as to guarantee that the Fourier kernel holds the same energy and admissibility coefficient of its generating wavelet.

\begin{proposition}
\label{proposition-3}
If $\psi(t)$ is a real wavelet and 
$\operatorname{\mathcal{H}}\{\psi\}(t)$ its Hilbert transform, then Ft{$\psi(t)$} is also a real wavelet with same energy and admissibility coefficient of its generating wavelet$\psi(t)$.
\end{proposition}
\proof
If $\psi(t)$ is a real wavelet and $\operatorname{\mathcal{H}}\{\psi\}(t)$ its Hilbert transform, then $\psi(t)$ and $\operatorname{\mathcal{H}}\{\psi\}(t)$ belong to~$L^2(\mathbb{R})$ and satisfy the admissibility condition. The energy and admissibility coefficient of $\operatorname{\mathscr{F}}\{\psi(t)\}$ are given by:
\begin{align*}
\operatorname{\mathcal{E}}
\left[
\operatorname{\mathscr{F}}\{\psi\}(t)
\right]
&
=
\frac{1}{2}
\int_{-\infty}^\infty
\left|
\psi(t)
-
j
\cdot
\operatorname{\mathcal{H}}\{\psi\}(t)
\right|^2
\operatorname{d}t
,
\\
\operatorname{C}
[
\operatorname{\mathscr{F}}\{\psi\}(t)
]
&=
\frac{1}{2}
\int_{-\infty}^\infty
\frac{\left|
\Psi(\omega)
-
\operatorname{sgn}(\omega)
\cdot
\Psi(\omega)
\right|^2}
{|\omega|}
\operatorname{d}\omega
.
\end{align*}
Simple manipulations furnish:
\begin{align*}
\operatorname{\mathcal{E}}
\left[
\operatorname{\mathscr{F}}\{\psi\}(t)
\right]
&
=
\frac{1}{2}
\int_{-\infty}^\infty
\left|
\psi(t)
\right|^2
+
\left|
\operatorname{\mathcal{H}}\{\psi\}(t)
\right|^2
\operatorname{d}t
,
\\
\operatorname{C}
[
\operatorname{\mathscr{F}}\{\psi\}(t)
]
&=
\begin{cases}
0,& \omega \geq0
\\
\int_{-\infty}^\infty
\frac{
2
\cdot
\left|
\Psi(\omega)
\right|^2}
{|\omega|}
\operatorname{d}\omega
,
&
\omega<0.
\end{cases}
\end{align*}

From Proposition~\ref{proposition-1}, 
it follows that 
$\operatorname{\mathscr{F}}\{\psi\}(t)$
also belongs to~$L^2(\mathbb{R})$ and it has the same energy as
$\psi(t)$. Since $\psi(t) \in L^2(\mathbb{R})$ and $|\Psi(0)|=0$,
it follows that $\operatorname{\mathscr{F}}\{\psi(t)\}$ also has the same admissibility coefficient of $\psi(t)$.
\endproof

The next proposition shows that Ft{$\psi(t)$} has same vanishing moments than its generating wavelet,$\psi(t)$.

\begin{proposition}
\label{proposition-4}
Let $\psi(t)$ be a wavelet with $N$ vanishing moments, then 
$\operatorname{\mathscr{F}}\{\psi\}(t)$
has [[-0also $N$ vanishing moments.
\end{proposition}
\proof
From Proposition 2, it follows that the $n$th moment of
$\operatorname{\mathscr{F}}\{\psi\}(t)$
is given by:
$$
M_n[\operatorname{\mathscr{F}}\{\psi\}(t)]
=
\left.
\frac{
\left[
\Psi(\omega)
-
\operatorname{sgn}(\omega)
\cdot
\Psi(\omega)
\right]^{(n)}}
{(\sqrt2)^n \cdot
(-2\cdot\pi\cdot j)^n}
\right|_{\omega=0}
,
$$
which can also be written as
$$
M_n[\operatorname{\mathscr{F}}\{\psi\}(t)]
=
M_n[\psi(t)]
+
\frac{1}{j}
\cdot
M_n\left[
\operatorname{\mathcal{H}}\{\psi\}(t)
\right]
.
$$
Then, 
$\operatorname{\mathscr{F}}\{\psi\}(t)$
has also $N$ null moments.
\endproof

In the frequency domain,
Fourier-Like wavelets are null for 
$\omega>0$.
For $\omega < 0$,
they have the magnitude response of the generating wavelet multiplied by a scalar factor.
This is a typical behavior of analytical signals.

\subsection{Analytic Wavelets}

Based on the results obtained from Fourier-Like wavelets it is actually simple to define analytic wavelets. An analytic function 
$\operatorname{\mathcal{A}}\{g\}(t)$
 is a complex signal designed by a real function $g(t)$ and its Hilbert transform $\operatorname{\mathcal{H}}\{g\}(t)$~\cite{ref4}. In the framework of wavelets, an analytic wavelet, 
$\operatorname{\mathcal{A}}\{\psi\}(t)$, can be defined by:
\begin{align*}
\operatorname{\mathcal{A}}\{\psi\}(t)
=
\frac{1}{\sqrt{2}}
\cdot
\left(
\psi(t)
+
j
\cdot
\operatorname{\mathcal{H}}\{\psi\}(t)
\right)
.
\end{align*}

Analytic wavelets have also same energy, admissibility coefficient and null moments than their generating wavelet,
$\psi(t)$. 
The proofs are similar to that of Propositions~\ref{proposition-3} and~\ref{proposition-4}.
In the frequency domain, analytic wavelets are null for $\omega<0$. For $\omega>0$, they have the magnitude response of the generating wavelet multiplied by a scalar factor.

\subsection{Computing the Wavelet Analysis of Asymmetrical Real Signals}

It may be expected that using Fourier-Like or analytic wavelets, even and odd parts of an asymmetric real signal can be better analyzed, respectively, by an even wavelet and its Hilbert transform, i.e. an odd wavelet.
In both cases, it will be necessary to perform a complex wavelet analysis. It is also possible to analyze both symmetries of a real signal using a real wavelet. In that case, the Hartley kernel should be invoked.

\section{The Hartley Kernel On The Wavelet Analysis}

The Hartley transform kernel, or Hartley kernel, is defined by the ``cosine and sine'' function: 
$\operatorname{cas}(t) = \cos(t) + \sin(t)$.
Recalling that
$\operatorname{\mathcal{H}}\{\sin\}(t)=\cos(t)$,
the Hartley kernel can also be written as
$\operatorname{cas}(t) = 
\cos(t) - \operatorname{\mathcal{H}}\{\cos\}(t)$
or
$\operatorname{cas}(t) = 
\sin(t) + \operatorname{\mathcal{H}}\{\sin\}(t)$.

This simple remark motivates the definition of Hartley-Like wavelets by taking the sum or the difference of a given real wavelet and its Hilbert transform. Let us define the Hartley kernel of a wavelet, or a Hartley-Like wavelet, 
$\operatorname{\mathscr{H}}\{\psi(t)\}$, by:
\begin{align}
\label{equation-10}
\operatorname{\mathscr{H}}\{\psi(t)\}
=
\frac{1}{\sqrt2}
\cdot
\left[
\psi(t)
\mp
\operatorname{\mathcal{H}}\{\psi\}(t)
\right]
.
\end{align}

The following proposition proves that 
$\operatorname{\mathscr{H}}\{\psi(t)\}$
is also a wavelet and that the factor  makes the Hartley kernel and its generating wavelet to have same energy and admissibility coefficient. Additionally, Proposition 6 shows that
$\operatorname{\mathscr{H}}\{\psi(t)\}$
has same null moments than their generating wavelet,$\psi(t)$.

\begin{proposition}
\label{proposition-5}
If $\psi(t)$ is a real wavelet and $\operatorname{\mathcal{H}}\{\psi(t)\}$ its Hilbert transform, then $\operatorname{\mathscr{H}}\{\psi(t)\}$ is also a real wavelet with same energy and admissibility coefficient of its generating wavelet,$\psi(t)$.
\end{proposition}
\proof
If $\psi(t)$ is a real wavelet and $\operatorname{\mathcal{H}}\{\psi(t)\}$ its Hilbert transform, then $\psi(t)$ and $\operatorname{\mathcal{H}}\{\psi(t)\}$ belong to~$L^2(\mathbb{R})$ and hold the admissibility condition. The energy and admissibility coefficient of $\operatorname{\mathscr{H}}\{\psi(t)\}$ are given by:

\begin{align*}
\operatorname{\mathcal{E}}
\left[
\operatorname{\mathscr{H}}\{\psi\}(t)
\right]
&
=
\frac{1}{2}
\int_{-\infty}^\infty
\left|
\psi(t)
\pm
\operatorname{\mathcal{H}}\{\psi\}(t)
\right|^2
\operatorname{d}t
,
\\
\operatorname{C}
[
\operatorname{\mathscr{H}}\{\psi\}(t)
]
&=
\frac{1}{2}
\int_{-\infty}^\infty
\frac{\left|
\Psi(\omega)
\pm
j
\cdot
\operatorname{sgn}(\omega)
\cdot
\Psi(\omega)
\right|^2}
{|\omega|}
\operatorname{d}\omega
.
\end{align*}
A simple manipulation yields to:
\begin{align*}
\operatorname{\mathcal{E}}
\left[
\operatorname{\mathscr{H}}\{\psi\}(t)
\right]
&
=
\frac{1}{2}
\int_{-\infty}^\infty
\left|
\psi(t)
\right|^2
+
\left|
\operatorname{\mathcal{H}}\{\psi\}(t)
\right|^2
\operatorname{d}t
,
\\
\operatorname{C}
[
\operatorname{\mathscr{F}}\{\psi\}(t)
]
&=
\int_{-\infty}^\infty
\frac{
\left|
\Psi(\omega)
\right|^2}
{|\omega|}
\operatorname{d}\omega
.
\end{align*}
From Proposition~\ref{proposition-1},
it follows that
$\operatorname{\mathscr{H}}\{\psi(t)\}$ also stays on 
$L^2(\mathbb{R})$ and it has same energy of $\psi(t)$. Additionally, $\operatorname{\mathscr{H}}\{\psi(t)\}$ has also same admissibility coefficient than its generating wavelet,
$\psi(t)$.
\endproof

\begin{proposition}
\label{propositon-6}
Let $\psi(t)$ be a wavelet with $N$ vanishing moments,
then $\operatorname{\mathscr{H}}\{\psi(t)\}$ has also $N$ vanishing moments.
\end{proposition} 
\proof
From Proposition~\ref{proposition-2},
it follows that the nth moment of $\operatorname{\mathscr{H}}\{\psi(t)\}$ is given by:
$$
M_n[\operatorname{\mathscr{H}}\{\psi\}(t)]
=
\left.
\frac{
\left[
\Psi(\omega)
\pm j
\cdot
\operatorname{sgn}(\omega)
\cdot
\Psi(\omega)
\right]^{(n)}}
{(\sqrt2)^n \cdot
(-2\cdot\pi\cdot j)^n}
\right|_{\omega=0}
,
$$
which can also be written as
$$
M_n[\operatorname{\mathscr{H}}\{\psi\}(t)]
=
M_n[\psi(t)]
\pm
M_n\left[
\operatorname{\mathcal{H}}\{\psi\}(t)
\right]
.
$$
Then, $\operatorname{\mathscr{H}}\{\psi(t)\}$ has also
$N$ null moments.

In the frequency domain, Hartley-Like wavelets have the magnitude response of the generating wavelet multiplied by a scalar factor. Additionally, they impose a 
$\pm\pi/4$-shift on the phase response of the generating wavelet.

\begin{figure*}
\centering

\subfigure[]{\includegraphics[width=0.3\linewidth]{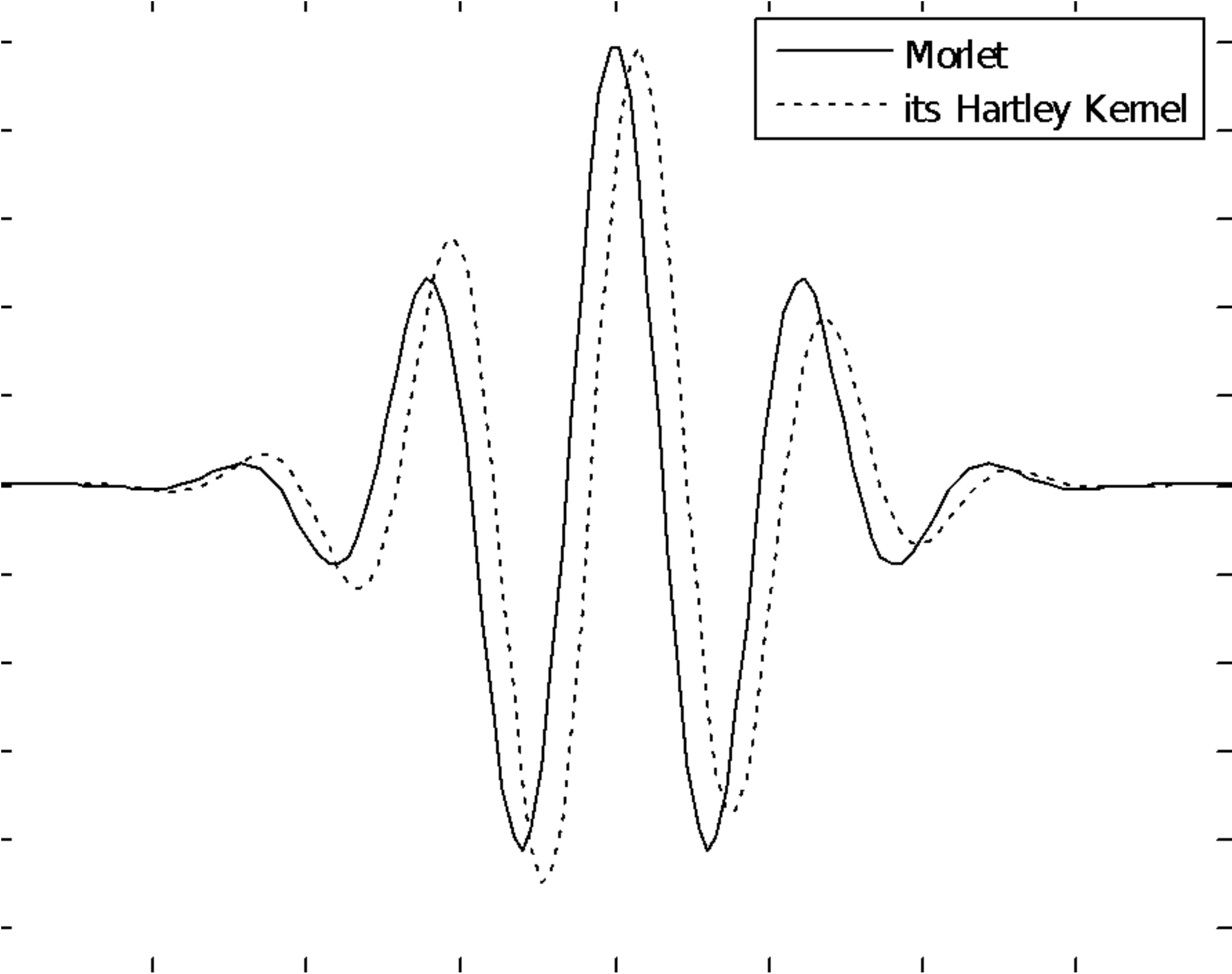}}
\quad
\subfigure[]{\includegraphics[width=0.3\linewidth]{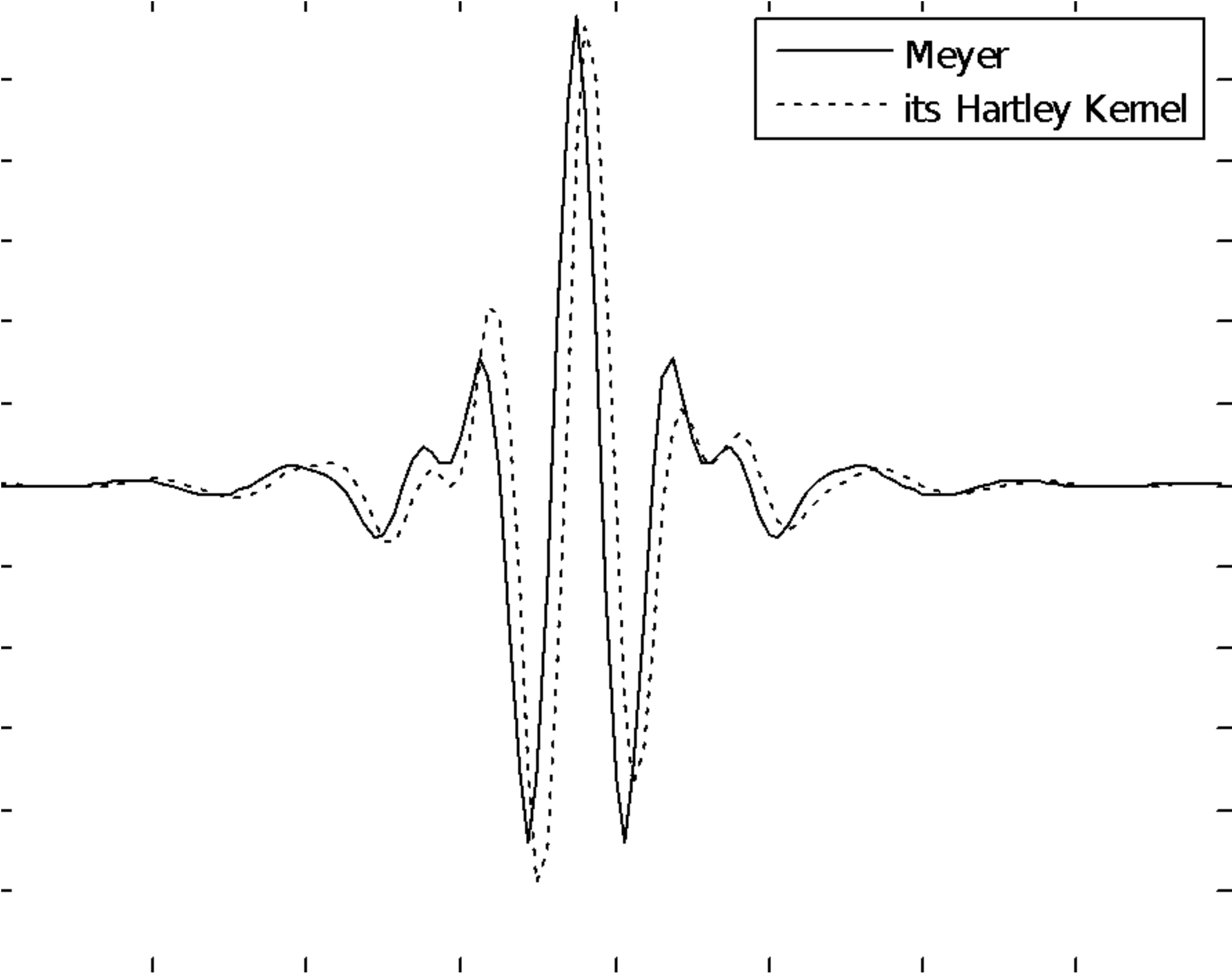}}
\quad
\subfigure[]{\includegraphics[width=0.3\linewidth]{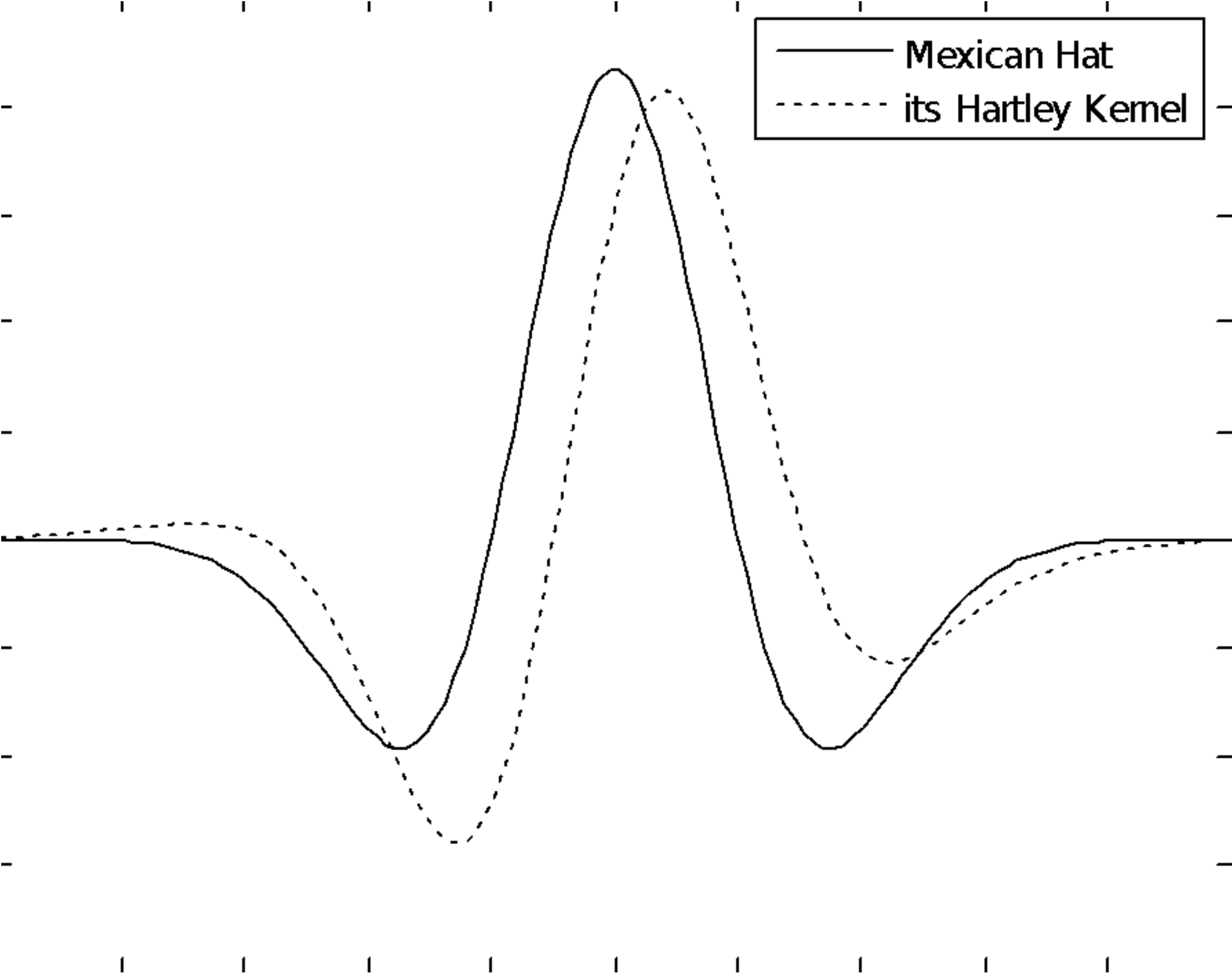}}
\\
\subfigure[]{\includegraphics[width=0.3\linewidth]{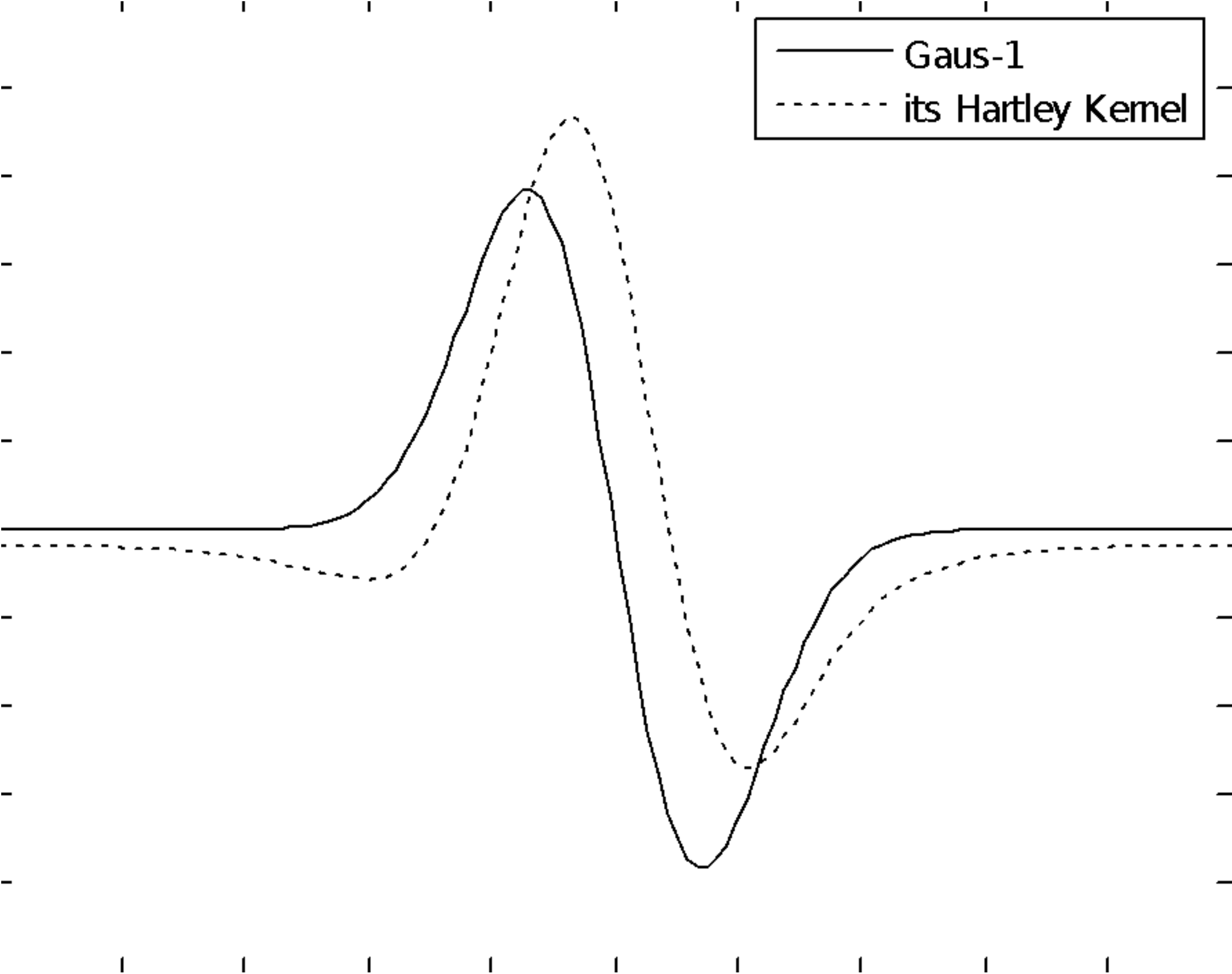}}
\quad
\subfigure[]{\includegraphics[width=0.3\linewidth]{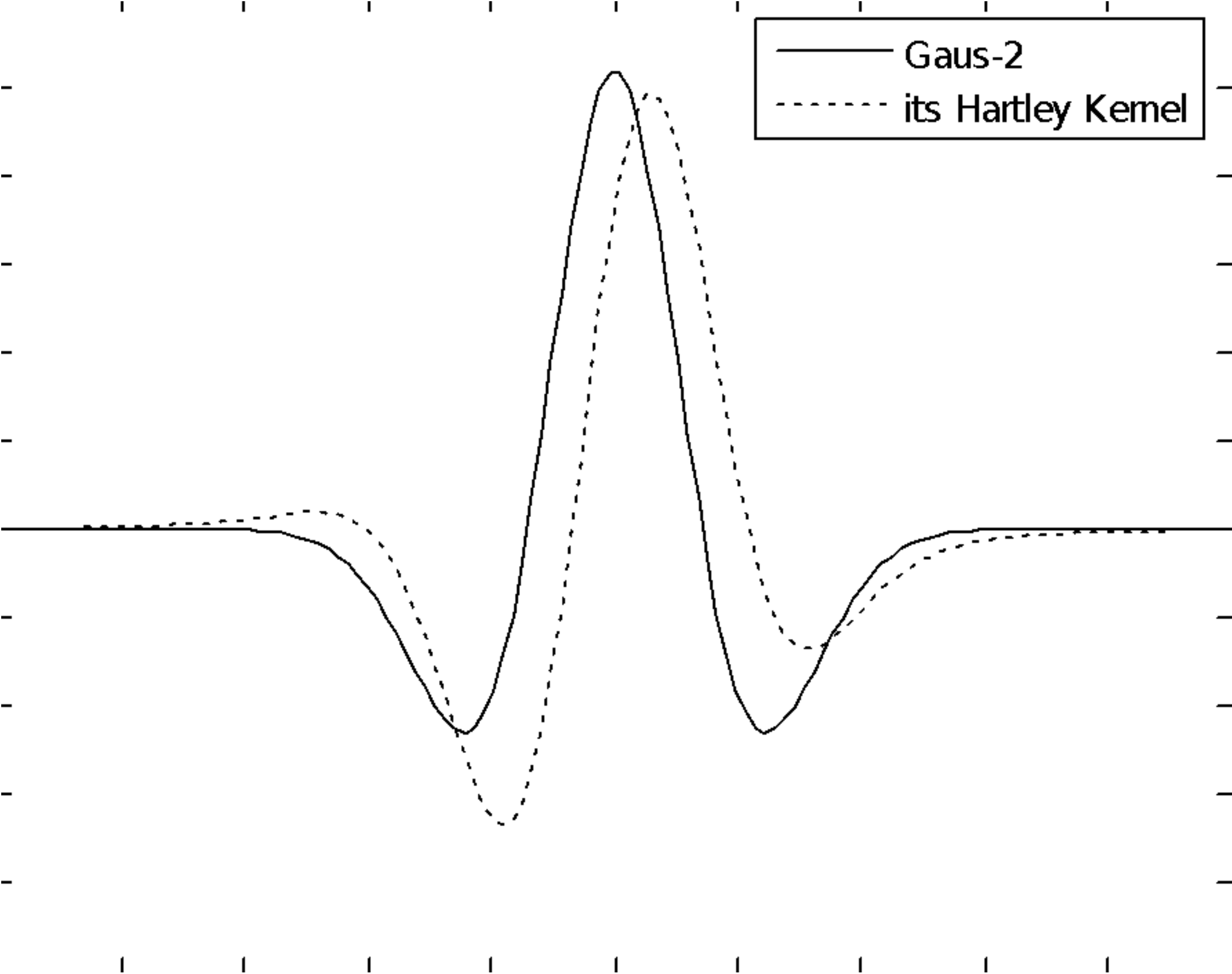}}
\quad
\subfigure[]{\includegraphics[width=0.3\linewidth]{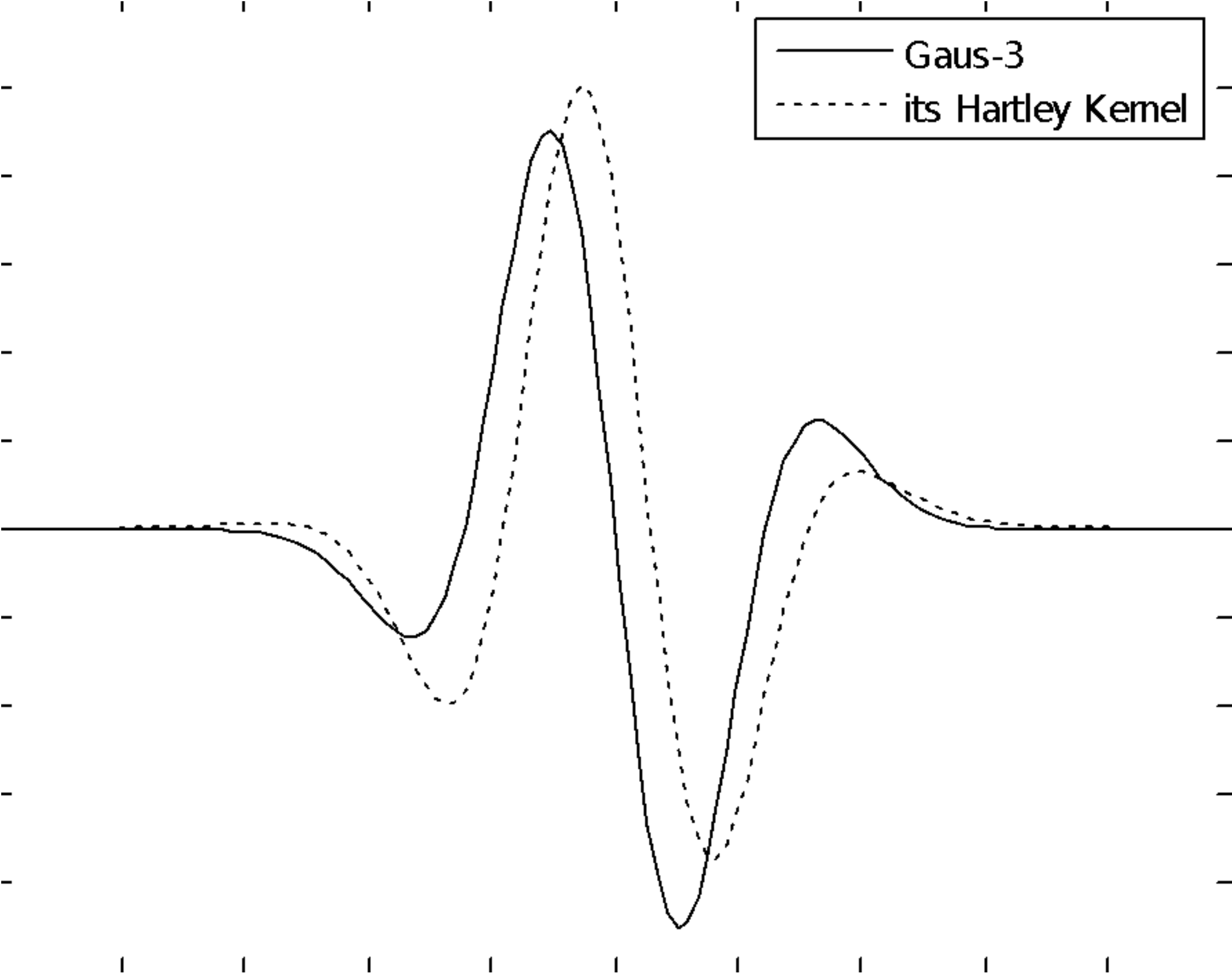}}

\caption{Continuous-time wavelets and their Hartley kernels: 
(a)~Morlet;
(b)~Meyer;
(c)~Mexican Hat;
(d)~Gaussian-1;
(e)~Gaussian-2;
and
(f)~Gaussian-3.
Scale of the horizontal axis is
one~unit (time) per division.
Scale of the vertical axis is
$0.2$ units (amplitude) per division,
ranging from $-1$ to~$1$.}
\label{figure-2}
\end{figure*}

Figure~\ref{figure-2}
shows some continuous-time wavelets and their corresponding Hartley kernels using 
the addition operator in~\eqref{equation-10}.

\section{Some Examples of Signal Analysis Using Fourier-Like and Hartley-Like Wavelets}

The wavelets proposed in this paper were simulated using the MATLAB Wavelet Toolbox~\cite{ref2}. Standard sample signals were analyzed to illustrate the behavior of the proposed wavelets.
Consider the wavelet transform  coefficients $C_{a,b}$ given by:
\begin{align*}
C_{a,b}
=
\frac{1}{\sqrt{a}}
\cdot
\int_{-\infty}^\infty
f(t)
\cdot
\psi
\left(
\frac{t-b}{a}
\right)
\operatorname{d}t
,
\end{align*}
where $a>0$ and $b$ are
real scale and translation scalars,
respectively.

\subsection{Applying the Hartley-Like Wavelet Analysis}

\begin{figure}
\centering

\includegraphics[width=0.7\linewidth]{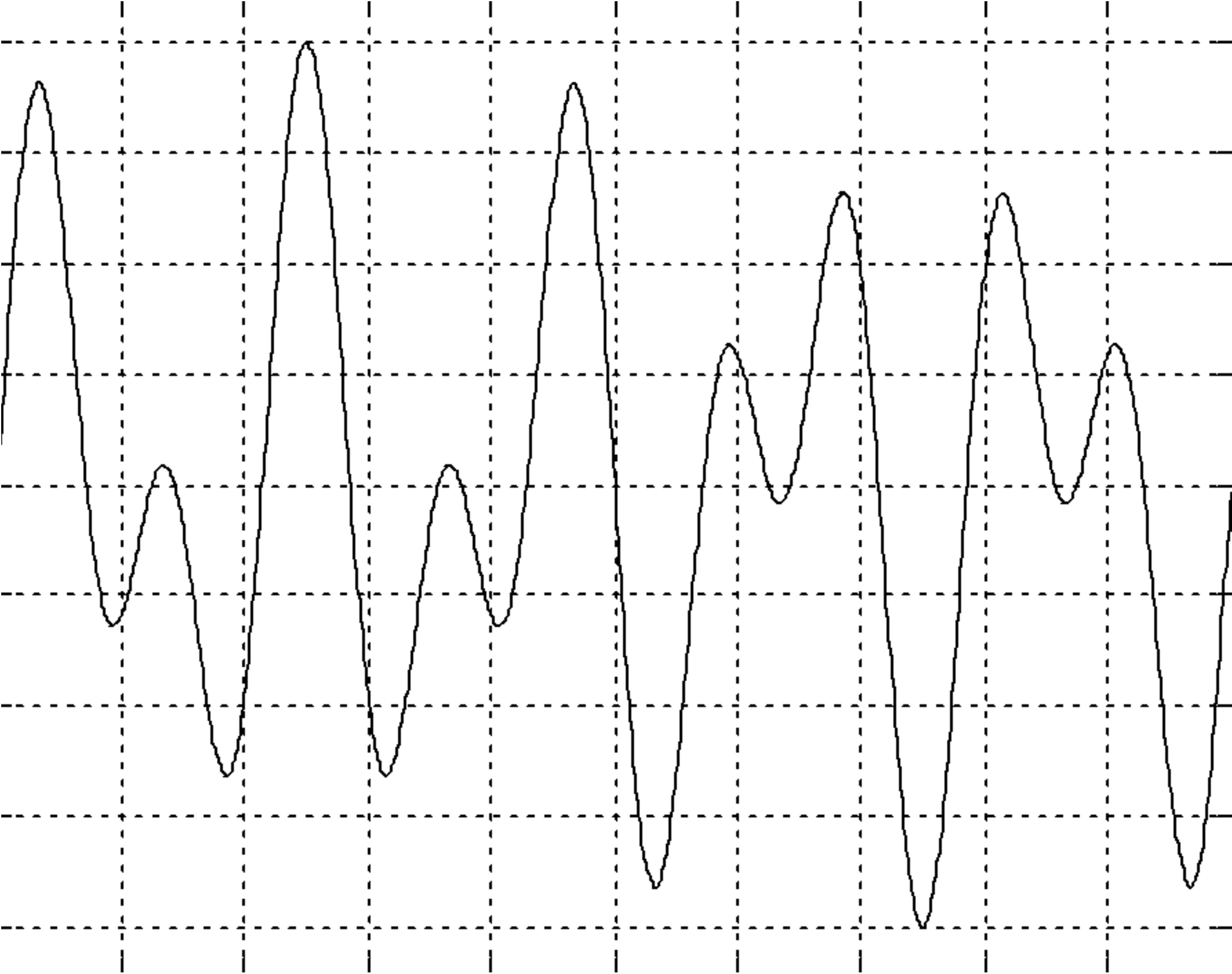}

\caption{A signal composed by two sinusoidal functions.
Scale of the horizontal axis is
100~unit (time) per division,
ranging from 0 to 1000.
Scale of the vertical axis is
$0.5$ units (amplitude) per division,
ranging from $-2$ to~$2$.}
\label{figure-3}
\end{figure}

\begin{figure*}
\centering

\subfigure[]{\includegraphics[width=0.3\linewidth]{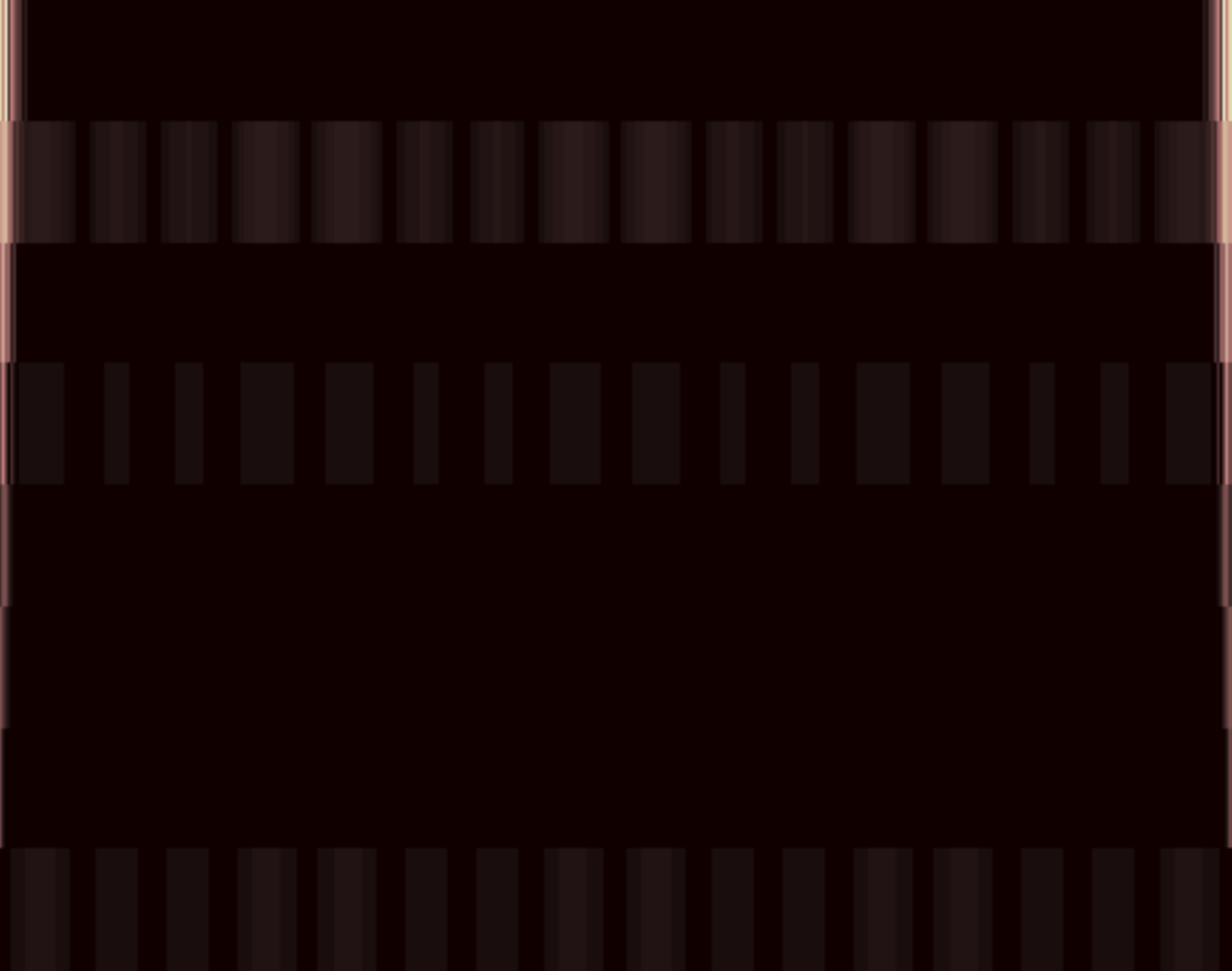}}
\quad
\subfigure[]{\includegraphics[width=0.3\linewidth]{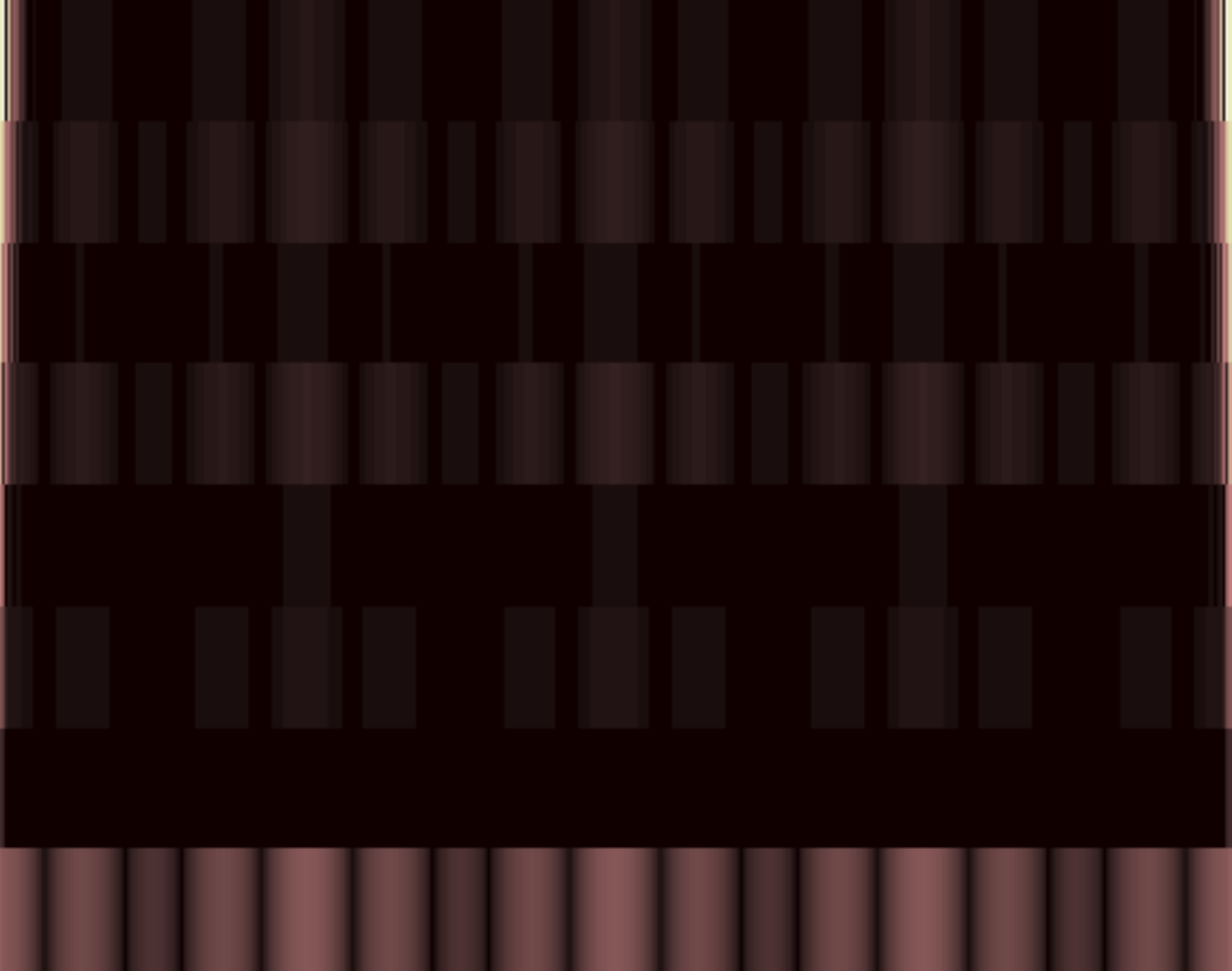}}
\quad
\subfigure[]{\includegraphics[width=0.3\linewidth]{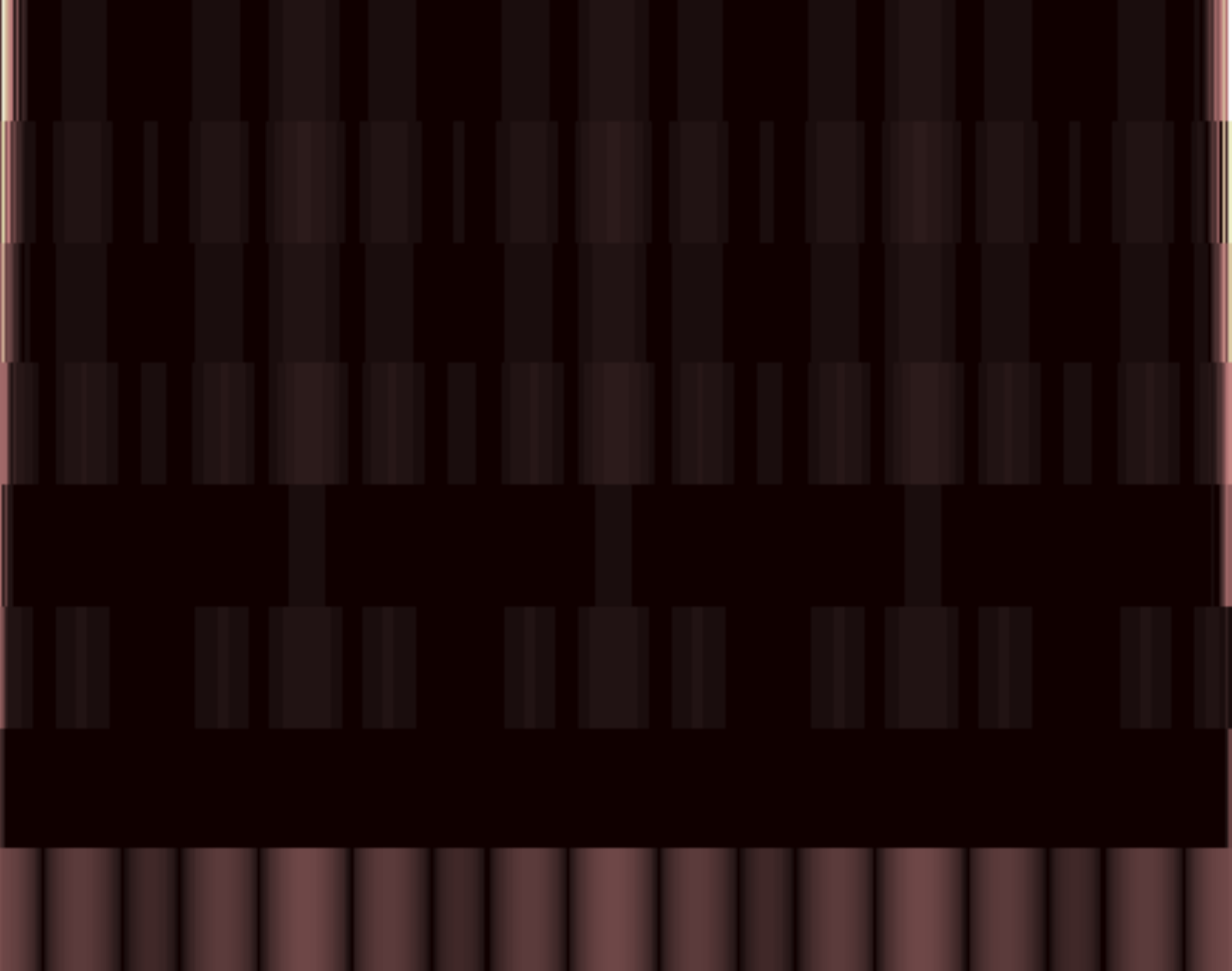}}

\caption{An 8-level wavelet analysis of later signal using: 
(a)~the Morlet wavelet; (b)~its Hilbert transform, and (c)~its Hartley kernel.
Scale of the horizontal axis is
100~unit (time or space~$b$) per division,
ranging from 0 to 1000.
Scale of the vertical axis is
one units (scale~$a$) per division,
ranging from $1$ to~$8$.}
\label{figure-4}
\end{figure*}

Figure~\ref{figure-3} shows a signal composed by two unitary-amplitude senoidal functions, 5~Hz and 9~Hz, and
Figure~\ref{figure-4} shows an 8-level wavelet analysis
(wavelet transform) of that signal using the Morlet wavelet, its Hilbert transform, and its Hartley kernel.

The scales where Hartley-Like Morlet wavelet may perform better signal analysis than the simple Morlet wavelet can be observed through the scales versus time charts.

From Figure~\ref{figure-4}, significant differences on the wavelet transform at level 1 can be seen, when using an even, an odd or an asymmetrical wavelet. It follows that, despite the time shifting impelled by the wavelet transform, asymmetrical wavelets can retrieve more information from a signal regardless its symmetry. As the kind of symmetry of the analyzing signal is not a priori known, the use of Hartley-Like wavelets on the continuous-time wavelet transform may achieve better results.

\begin{figure}
\centering

\includegraphics[width=0.7\linewidth]{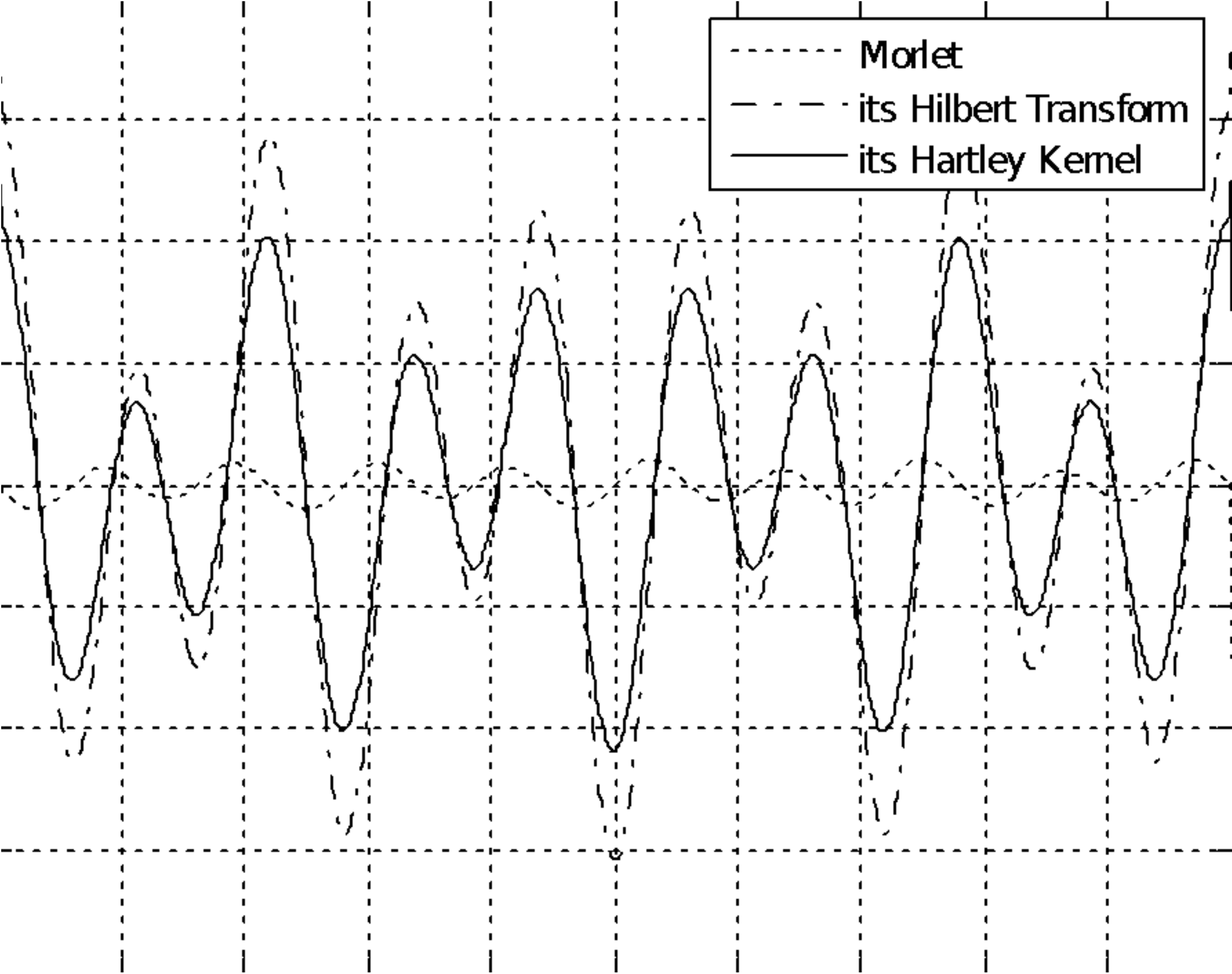}

\caption{The 1st-level wavelet transform of a combined signal using the Morlet wavelet, its Hilbert transform and its Hartley kernel.
Scale of the horizontal axis is
100~unit (time) per division,
ranging from 0 to 1000.
Scale of the vertical axis is
$0.005$ units (amplitude) per division,
ranging from $-0.02$ to~$0.02$.}
\label{figure-5}
\end{figure}

Figure~\ref{figure-5} shows the signal derived from the 1st-level wavelet transform of the later signal using this approach. Some improvements can be achieved using an odd or an asymmetrical wavelet to analyze an odd signal, than using an even wavelet.

\subsection{Applying the Analytic (or Fourier-Like) Wavelet Analysis}

\begin{figure}
\centering

\includegraphics[width=0.7\linewidth]{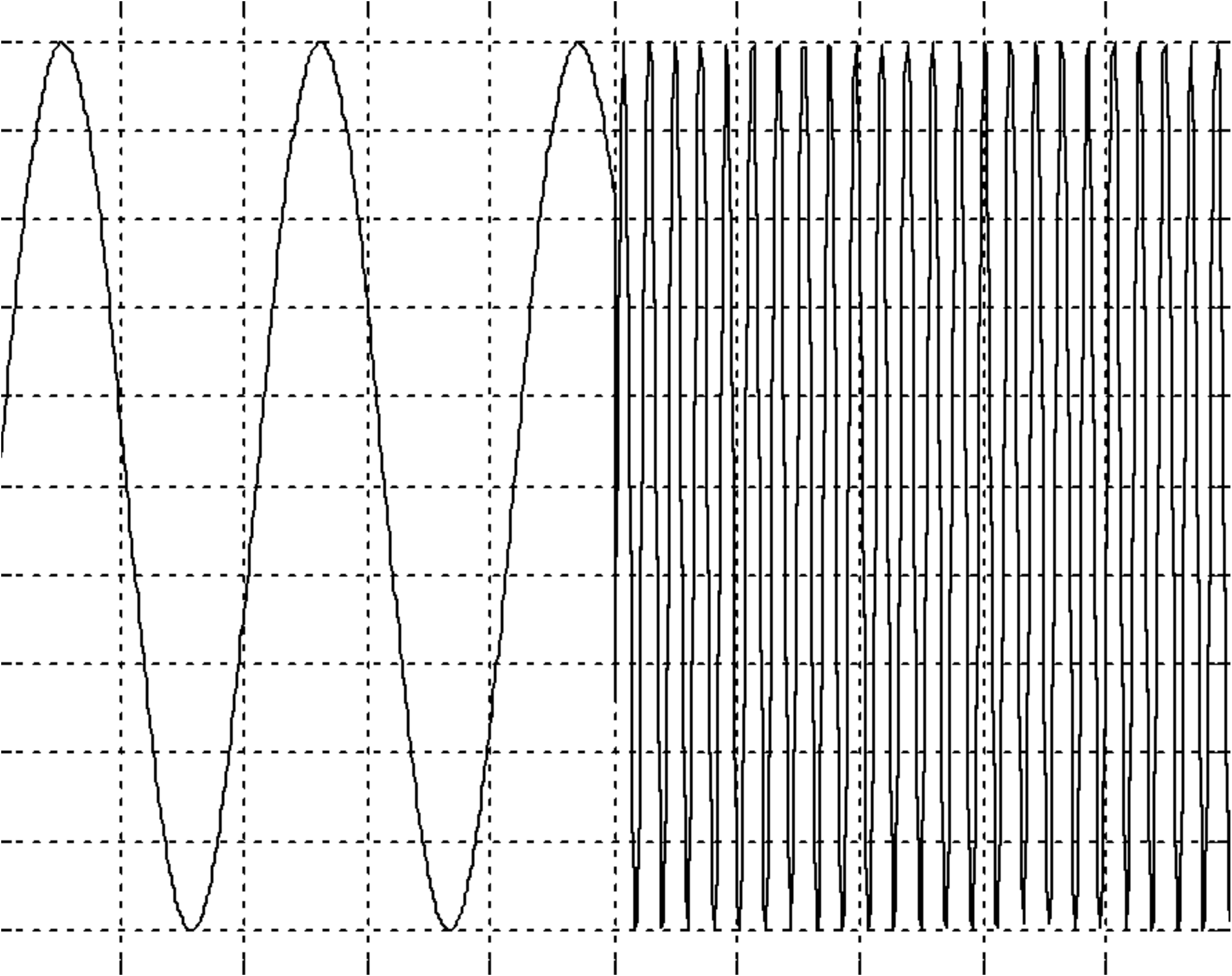}

\caption{A standard frequency breakdown signal.
Scale of the horizontal axis is
100~unit (time) per division,
ranging from 0 to 1000.
Scale of the vertical axis is
$0.2$ units (amplitude) per division,
ranging from $-1$ to~$1$.}
\label{figure-6}
\end{figure}

\begin{figure*}
\centering

\subfigure[]{\includegraphics[width=0.22\linewidth]{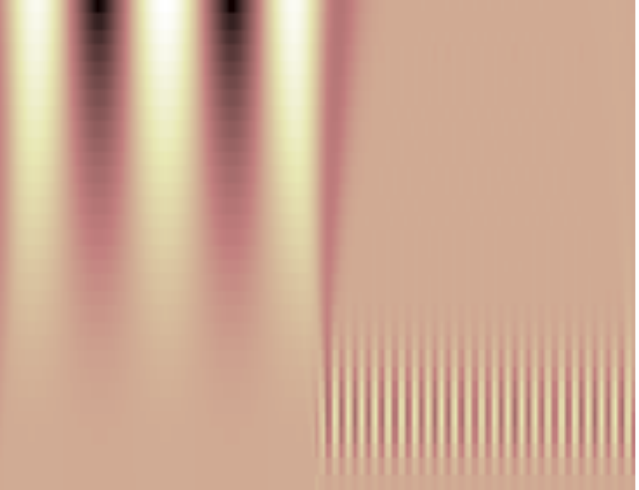}}
\quad
\subfigure[]{\includegraphics[width=0.22\linewidth]{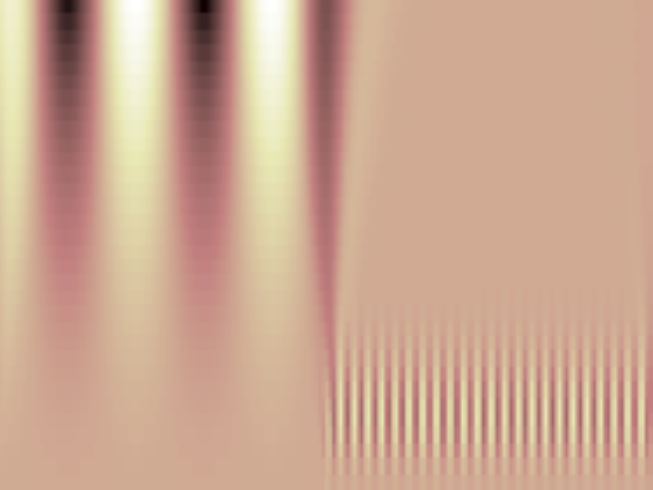}}
\quad
\subfigure[]{\includegraphics[width=0.22\linewidth]{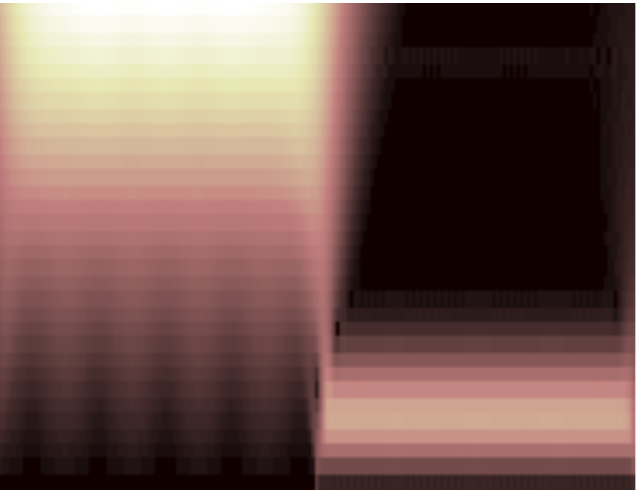}}
\quad
\subfigure[]{\includegraphics[width=0.22\linewidth]{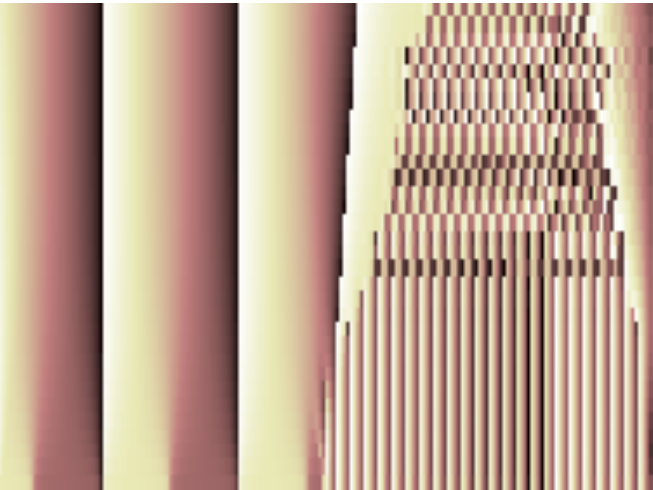}}
\caption{A 32-level wavelet analysis of a standard frequency breakdown signal using the Analytic Mexican Hat wavelet:
(a)~real part of the wavelet coefficients;
(b)~imaginary part of the wavelet coefficients;
(c)~magnitude of the wavelets coefficients;
and
(d)~phase of the wavelet coefficients.
Scale of the horizontal axis is
100~units (time or scale~$b$) per division,
ranging from 0 to 100.
Scale of the vertical axis is
one unit (scale) per division,
ranging from 1 to 31.}
\label{figure-7}
\end{figure*}

Hilbert transform has already been utilized in wavelet-based signal analysis~\cite{ref5}. In this section we carry a preliminary investigation about analytic wavelet analysis. Figure~\ref{figure-6} shows a standard frequency breakdown signal and Figure~\ref{figure-7} presents a 32-level wavelet analysis, assessed by the analytic Mexican Hat wavelet.

From Figure~\ref{figure-7}, it can be observed that the Mexican Hat wavelet can identify the presence of both frequencies as well the time when occurs the frequency changing, which can be better accurate at lowest scale values. The modulus of that analytic wavelet analysis shows that the high frequency signal can be viewed at lowest scale values and the low frequency signal at higher scale values.

\begin{figure}
\centering

\includegraphics[width=0.7\linewidth]{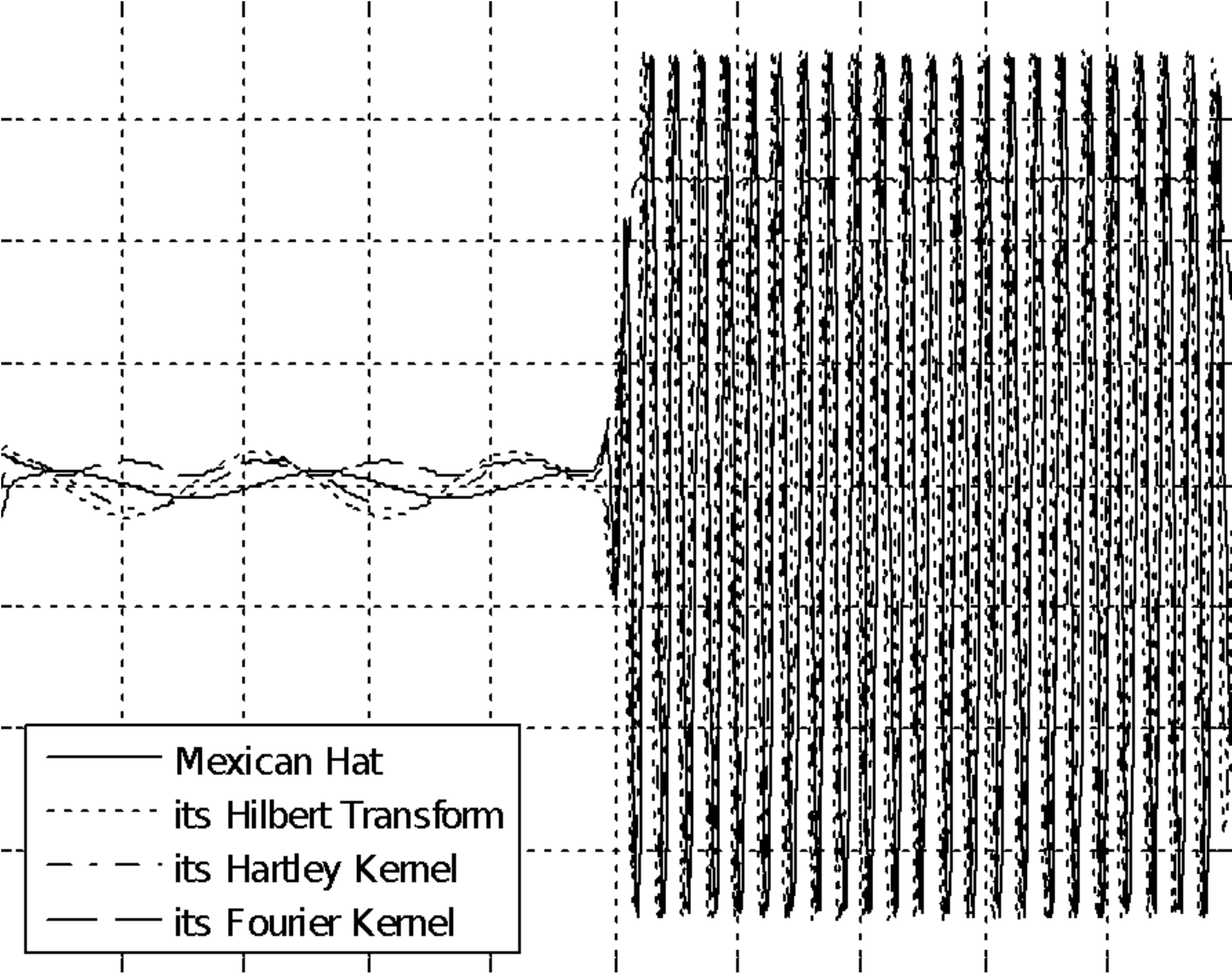}

\caption{The continuous-time wavelet transform of a standard frequency breakdown signal using the Mexican Hat wavelet, its Hilbert transform, its Hartley kernel, and its Fourier kernel as mother wavelets and assuming 5 as the scale parameter. The modulus is used for the Fourier-Like wavelet analysis.
Scale of the horizontal axis is
100~units (time) per division,
ranging from 0 to 1000.
Scale of the vertical axis is
one unit (amplitude) per division,
ranging from $-4$ to~$4$.}
\label{figure-8}
\end{figure}

\begin{figure}
\centering

\includegraphics[width=0.7\linewidth]{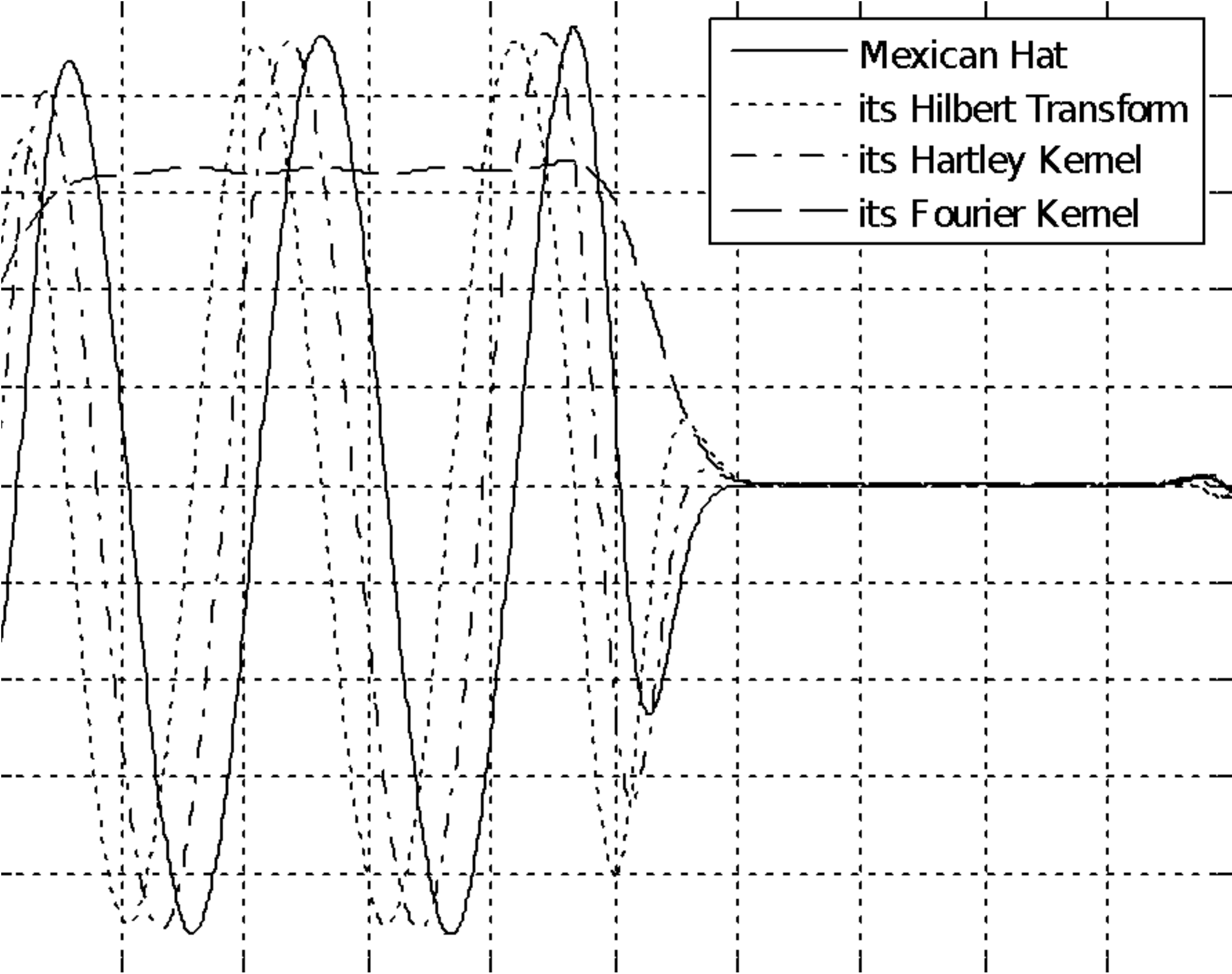}

\caption{The continuous-time wavelet transform of a standard frequency breakdown signal using the Mexican Hat wavelet, its Hilbert transform, its Hartley kernel, and its Fourier kernel as mother wavelet and assuming 25 as the scale parameter. The modulus is used for the Fourier-Like wavelet analysis.
Scale of the horizontal axis is
100~units (time) per division,
ranging from 0 to 1000.
Scale of the vertical axis is
one unit (amplitude) per division,
ranging from $-5$ to~$5$.}
\label{figure-9}
\end{figure}

As an example, Figures~\ref{figure-8} and~\ref{figure-9} show, respectively, the 5th and 25th-level wavelet transform, considering wavelets derived from the Mexican Hat wavelet. It can be seen the low and high frequency signals, when using real wavelets. 

\begin{figure}
\centering

\includegraphics[width=0.7\linewidth]{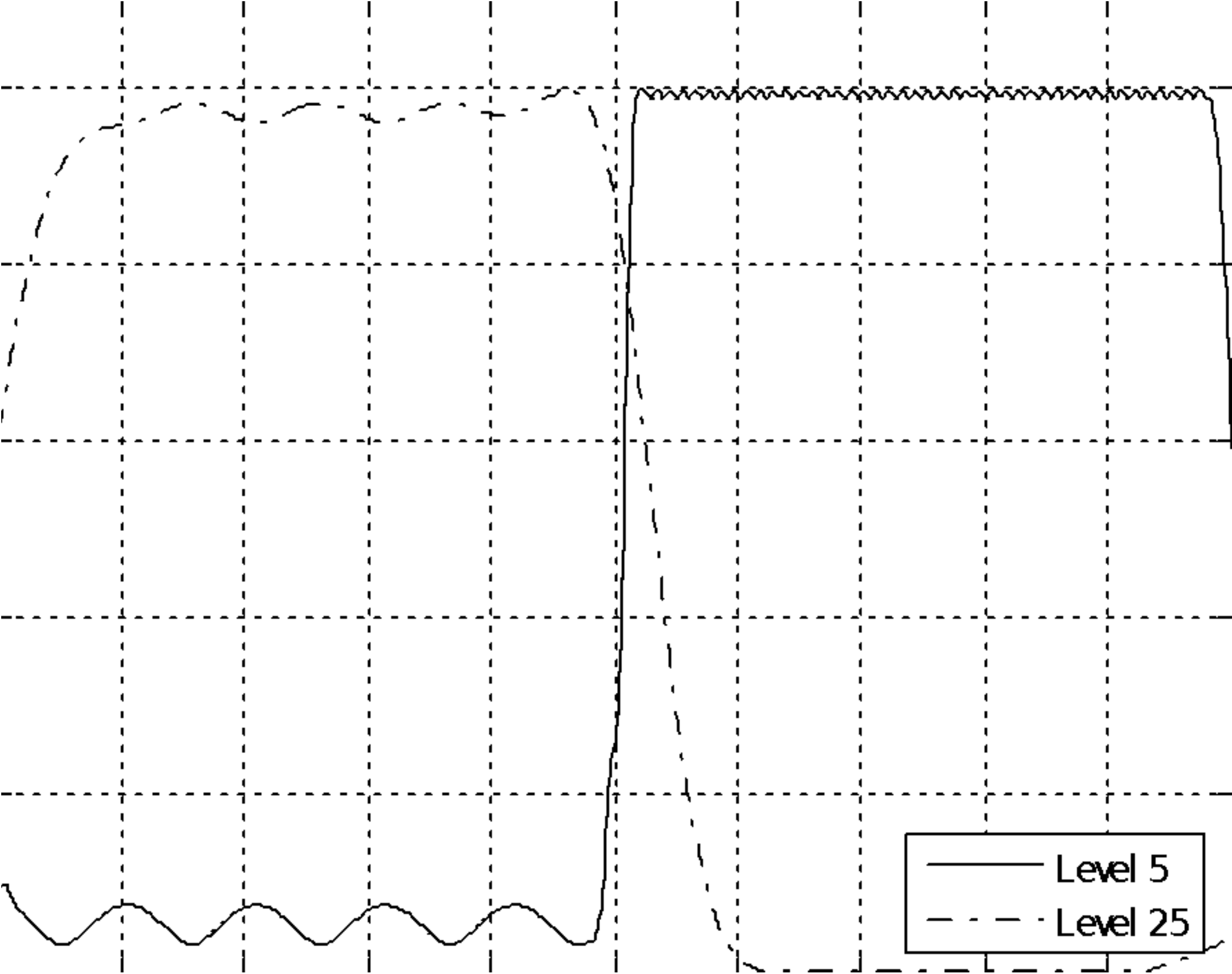}

\caption{Normalized modulus of the 5th and 25th-level wavelet coefficients of a standard frequency breakdown when using the analytic, or Fourier-Like, Mexican Hat wavelet.
Scale of the horizontal axis is
100~units (time) per division,
ranging from 0 to 1000.
Scale of the vertical axis is
$0.2$ unit (amplitude) per division,
ranging from $0$ to~$1$.}
\label{figure-10}
\end{figure}

Figure~\ref{figure-10} shows the normalized modulus of the 5th and 25th-level wavelet analysis, when using the analytic, or Fourier-Like, Mexican Hat wavelet. This feature shows that when analyzing individually scales, the time interval when occurs different frequency contents can be estimated.

\section{Conclusions}

New wavelet functions have been introduced that might improve the continuous-time signal analysis. Anti-symmetrical wavelets can be designed throughout the Hilbert transform of a symmetrical wavelet, and vice-versa. Together, they are called a Hilbert transform pair of wavelets. Fourier-Like and Hartley-Like wavelets have been derived from the Fourier and Hartley transform kernels, which have been written, on the basis of Hilbert transform.

The example cases illustrate that, despite the time shifting impelled by the wavelet transform, asymmetrical wavelets can grasp more features from a signal having no particular symmetry. Additionally, Fourier-Like and analytic wavelets have potential applications for disturbances and frequency detection.

{\small
\bibliographystyle{IEEEtran}
\bibliography{ref-clean}
}

\end{document}